\documentclass[9pt]{amsart}
\usepackage{amsmath}
\usepackage{amsthm}
\usepackage{epsfig}
\usepackage{amsfonts}
\usepackage {epsfig}

\swapnumbers \theoremstyle{plain}

\newtheorem{thm}{Theorem}[section]

\newtheorem{lem}[thm]{Lemma}

\newtheorem{prop}[thm]{Proposition}

\theoremstyle{definition}
\newtheorem{defn}[thm]{Definition}
\newtheorem{construction}[thm]{Construction}

\newtheorem{claim}[thm]{Claim}

\theoremstyle{remark}
\newtheorem*{rem}{Remark}

\swapnumbers

\newcommand{\R}{$\mathbf{R}P^3 $}
\newcommand{\GG}{\ensuremath{\overline{G_{1, e_{1}}}}}
\newcommand{\GGG}{\ensuremath{\overline{G_{2, e_{2}}}}}

\title{ Connected Sums of Closed Orientable Triangulated 3-Manifolds}
\author{Alexander Barchechat}
\date{\today}
\thanks{The author was partially supported by NSF grants XXXXX and by UC Davis}

\address{mathematics department, tulane university, new orleans, la, 70118}
\email{alexb@math.tulane.edu}

\usepackage{amssymb}

\begin{document}

\begin{abstract} In this paper, we describe geometrical constructions to obtain triangulations of connected
sums of closed orientable  triangulated 3-manifolds. Using these constructions, we show that it takes time
 polynomial in the number of tetrahedra to check if a closed orientable 3-manifold, equipped with a minimal
 triangulation, is reducible or not. This result can easily be generalized to compact orientable 3-manifolds with
 non-empty boundary.
\end{abstract}

\maketitle

\section{Introduction}
\subsection{Normal Surfaces}

This paper relies heavily on the process of ``collapsing'' normal 2-spheres, which was first described by Jaco
and Rubinstein in \cite{JR}. In section 2, we describe an alternate ``collapsing'' process, based on the existence
of {\it collapsing surfaces}. In section 3, we give a proof of the main theorem:

Theorem 2.1:
 {\it Let $M$ be a closed orientable triangulated 3-manifold with $t$ tetrahedra. Let S be a non-trivial
  normal 2-sphere. Then M  is homeomorphic to $ M_{1} \# M_{2}$...$\# M_{k}$ $\# r_{1}(S^{1}\times S^{2})$
   $\# r_{2}\mathbf {RP}^{3}$ $\#r_{3}L(3,1)$, where $r_{1}$, $r_{2}$, $r_{3}$, $k$ $\geq 0$,
   $|M_{1}| +...+ |M_{k}| < |M|$ and the $M_{i}$'s are closed orientable triangulated 3-manifolds.}

 In section 4, we describe
geometrical constructions to obtain triangulations for connected
 sums of closed orientable triangulated 3-manifolds. In section 5, we prove our main result regarding the
 existence of normal 2-spheres with at most 2 non-zero quadrilateral types in minimal triangulations. Finally,
 in the last section, we fully describe an algorithm, due to Andrew Casson,
  to decompose a closed orientable triangulated 3-manifold into irreducible pieces. Using this algorithm, we show
   how one can check, in time polynomial with respect to the number of tetrahedra in the triangulation, if
  a closed orientable 3-manifold equipped with a minimal triangulation is reducible or not.

All the 3-manifolds considered in this paper are in the piecewise linear category, i.e. every 3-manifold will be
associated with a triangulation.

\begin{defn} A \textbf{\textit{pseudo-triangulation}} of a compact orientable 3-manifold $M$ is a set $\Delta$ of pairwise
disjoint tetrahedra, together with a family
of homeomorphisms, $\Phi$, where the domain and image of each homeomorphism consist of faces of tetrahedra. The
identification space, $\Delta$/$\Phi$,
is homeomorphic to $M$. If $\Delta$ consist of n tetrahedra,
 we write $\mathbf {|M|} = n$. A \textbf{\textit{genuine triangulation}} is a pseudo-triangulation such that every
 tetrahedron is embedded in M. By abuse of language, we call a pseudo-triangulation a triangulation. Moreover, M
 will always be associated with a fixed triangulation.
\end{defn}

We denote the i-skeleton of M by $T^{(i)}$, for $i=0$, 1, 2, 3.
  Let $R: M \times I \rightarrow M$ be an isotopy of $M$. $R$ is called a
  \textbf{\textit{normal isotopy}} if it is invariant in each tetrahedron $\Delta_i$, i.e. $R$
  ($\Delta_i$, t)= $\Delta_i$ for all $t$ $\in I$.  \textbf{\textit{Normal surfaces}}, which were first developped
  by W.Haken ~\cite{Ha},
are embedded surfaces in $M$ which intersect each tetrahedron in planes in general position with
 respect to $T^{(2)}$. There are 7 different isotopy classes of planes (called \textbf{\textit{elementary disks}}
 ) for each
  tetrahedron: 4 triangle types and 3 quadrilateral types. A normal surface can thus be described as an ordered set of elementary
   disks in each tetrahedron. Using this definition, a normal surface describes a whole class of
   embedded surfaces which are all normal isotopic. Hence, by a normal surface we mean its normal isotopy class.

We call a normal 2-sphere \textbf{\textit{trivial}} if it intersects the tetrahedra in triangles only, and we call it
 non-trivial otherwise.

Given a closed embedded normal surface S, it has to satisfy the  \textbf{\textit{quadrilateral property}}:
if there is a quadrilateral type in a tetrahedron, then no other types of quadrilateral can exist in that same tetrahedron.
Moreover, F must satisfy \textbf{\textit{matching equations}}:

If $M$ has $t$ tetrahedra, there are exactly $6t$ matchings to be satisfied, 3 for each face of the
 triangulation. We have noticed above that a normal surface is described by a set of $7t$ elementary
  disks, 7 for each tetrahedron. So one way to think of a normal surface is to look at it as a non-negative
   integer valued vector with $7t$ entries satisfying the quadrilateral property and the matching equations.
    Given a normal surface $F$, we may refer to it as $x_F$, where $x_F = (x_1, x_2,$... $,
    x_{7t})$. The matching equations may be seen as a set
    of $6t$ equations of the form $x_{i} + x_{j} = x_{k} + x_{l}$.


Let $F$ and $G$ be two normal surfaces which intersect each other. Let $x_{F}$ and $x_{G}$ be their corresponding
 vectors. Suppose that $F$ and $G$ intersect each other with the only requirement that their quadrilateral
 types are the same in each tetrahedron. We can perform an operation, called a \textbf{\textit{regular exchange}},
 described in one of
  many ways in figure~\ref{fig1.3} below, such that the resulting pieces are disjoint elementary disks.

\bigskip
\begin{figure}[h]
\hspace{1in} \psfig{figure=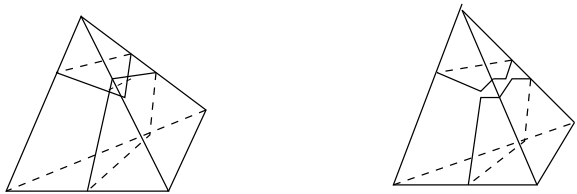, width= 3in} \caption{Before and after a regular exchange.} \label{fig1.3}
\end{figure}

The surface obtained in this manner is the normal surface represented
by the vector $x_{F} + x_{G}$.
 This operation is called the \textbf{\textit{surface addition}} or \textbf{\textit{Haken sum}} of $F$
and $G$, and it is unique up to isotopy.
 We define the  \textbf{\textit{weight}} of a normal surface $F$, denoted by \textbf{\textit{wt(F)}}, as the number
     of intersection points of $F$ with $T^{(1)}$. It is well known that the Haken sum
      preserves the weight and the Euler characteristic of normal surfaces.

A normal surface is called \textbf{\textit{fundamental}} if it cannot be written as the sum of two non-isotopically
parallel surfaces. Obviously, a fundamental surface has to be connected. Haken proved that the set of fundamental
surfaces is finite and that they can be found algorithmically (see ~\cite{Hak} and~\cite{Ha}).

Consider the set of non-negative real solutions to the matching equations and satisfying the quadrilateral property.
 It is well known, through linear programming, that this set forms a cone in $\mathbf{R}^{7t}$. We intersect this cone with
 the set of solutions to the equation: $\sum_{i=1}^{7t}x_i = 1$. We obtain a convex polyhedron called the
 \textbf{\textit{projective solution space}} of $M$ with respect to its triangulation $T$, and we denote
it by \textit{P(M, T)}. For each normal surface $S$, there corresponds a rational vector \=S in \textit{P(M, T)}
  called the projective class of $S$. Conversely, any rational vector in \textit{P(M, T)} can be multiplied by
   an integer to obtain a vector representing a normal surface.

It can be shown that the vertices of \textit{P(M, T)} have rational entries. Let $v$ be a vertex of
\textit{P(M, T)} and let $k$ be the smallest integer such that $k \cdot v$ is an integral solution. We call
$k \cdot v$ a \textbf{\textit{vertex solution}}. In particular, an integral solution $F$ is a vertex solution
 if and only if the integral solutions, $X$ and $Y$, to the equation $n \cdot F = X + Y$ are multiples of $F$.
  We call $F$ a \textbf{\textit{vertex surface}} if it is connected, 2-sided, and if its representative on
  \textit{P(M, T)} is a vertex. Note, if $F$ is a vertex surface, then either $F$ is also a vertex solution
  or it is the double of a vertex solution. See  ~\cite {JT},~\cite{Ha},
  ~\cite{HLP},  ~\cite{JR}, ~\cite{JS},
  and ~\cite{Mat} for more details on normal surfaces.

\section{Collapsing normal 2-spheres}
\label{collapse}

This work is directly inspired by a theorem, stated first by Jaco and Rubinstein, which appeared in ~\cite{JR}. The
original theorem is the following:

{\it Theorem (Jaco-Rubinstein)}: Let $M$ be a closed orientable triangulated 3-manifold. If $M$ contains a
non-trivial 2-sphere, then either $M \cong M_{1} \# M_{2}$, with $|M_{1}| + |M_{2}| < |M|$, or
$M \cong M_{1} \# r_{1}(S^{1} \times S^{2}) \# r_{2}\mathbf{RP}^{3} \# r_{3}L(3, 1)$, with $|M_{1}| < |M|$.

In this section, we define the necessary terminology to prove the following theorem.

\begin{thm} \label{MT}
 Let $M$ be a closed orientable triangulated 3-manifold with $t$ tetrahedra. Let S be a non-trivial
  normal 2-sphere. Then M  is homeomorphic to $ M_{1} \# M_{2}$...$\# M_{k}$ $\# r_{1}(S^{1}\times S^{2})$
   $\# r_{2}\mathbf {RP}^{3}$ $\#r_{3}L(3,1)$, where $r_{1}$, $r_{2}$, $r_{3}$, $k$ $\geq 0$,
   $|M_{1}| +...+ |M_{k}| < |M|$ and the $M_{i}$'s are closed orientable triangulated 3-manifolds.
\end{thm}

If we cut $M$ along $S$, we obtain a cell decomposition of a
3-manifold $M \backslash S$ with two 2-spheres as boundary. The idea is to collapse each of the two 2-spheres to a
point, and obtain a well-defined triangulation for the resulting
3-manifold $M \backslash S$. After cutting $M$ along $S$, there
are 7 different types of polyhedra in the cell decomposition of $M
\backslash S$ to consider: tetrahedra, truncated tetrahedra (with
1, 2, 3 or 4 truncations),
 prisms, truncated prisms (with 1 or 2 truncations), tips, $I \times$ quadrilateral, and $I \times$ triangles. The
  last two polyhedra are
  called \textbf{\textit{I-bundles}}. Some types of polyhedra may be combinatorially equivalent (e.g.
  a tetrahedron and a tip), but for topological reasons we consider them as different.

\bigskip
\begin{figure}[h]
\hspace{.3in} \psfig{figure=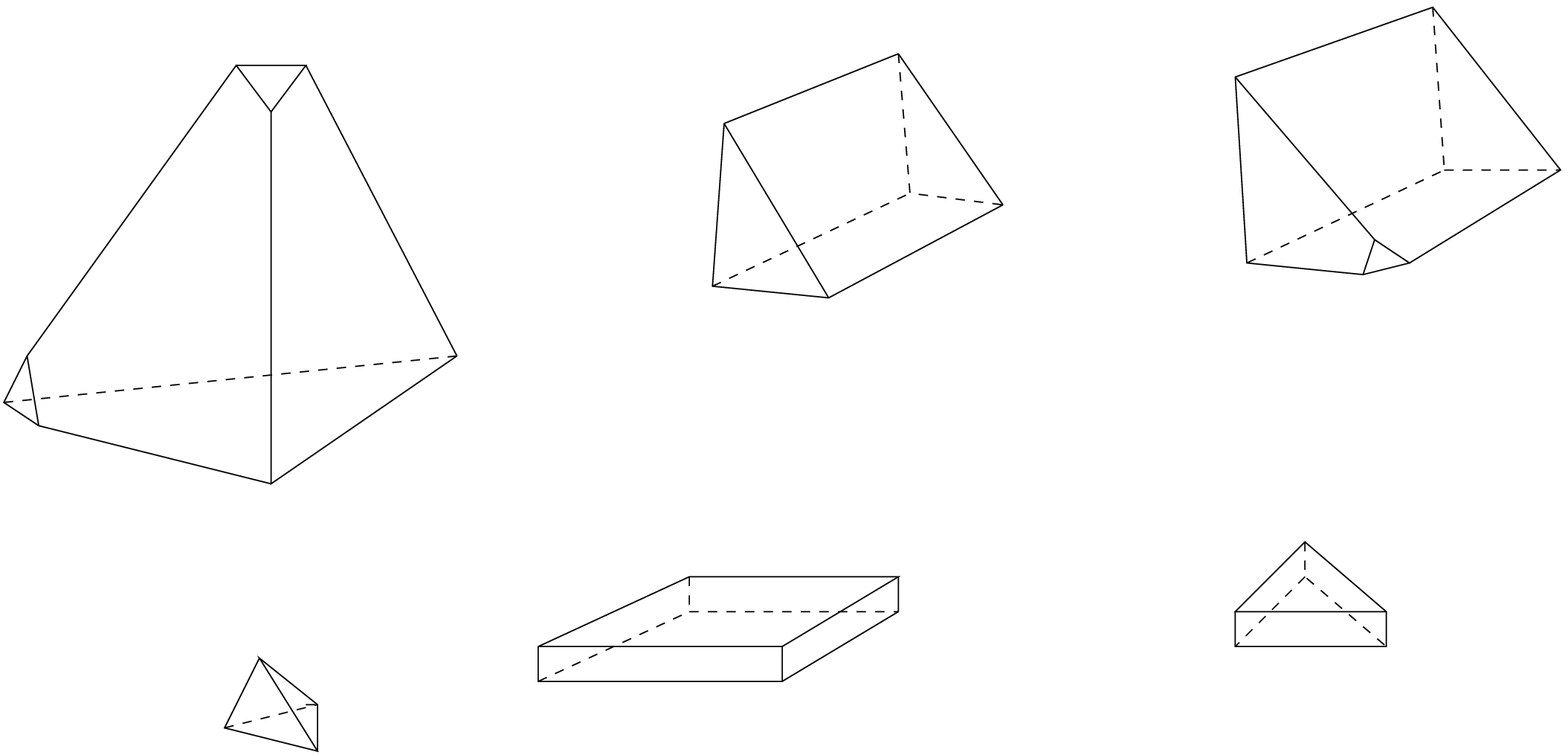, width= 3in} \caption{Six types of polyhedra in $M \backslash S$} \label{fig1}
\end{figure}
\bigskip

\begin{defn}
Let $M$ be a closed orientable triangulated 3-manifold and let $S$ be an embedded orientable surface.
 We define $\mathbf{M} \backslash \mathbf{S}$ to be the 3-manifold $\overline{M-Nbhd(S)}$,
 where $Nbhd(S)$ denotes a regular neighborhood of $S$. By \textbf{\textit{collapsing}} $S$ in $M$ we mean
 collapsing the two copies of $S$ in $M \backslash S$ to points. If $S$ is a 2-sphere, then
 collapsing $S$ is topologically equivalent to cutting $M$ along $S$
 and capping off with 3-balls.

From now on, $S$ denotes a non-trivial normal 2-sphere unless
otherwise stated. We now describe the process of  collapsing a
prism. To do so, we define
 a prism $P$ as the quotient space
  $(I \times J) \times K / (a, 1, c_{1}) \sim (a, 1, c_{2})$, where
  $I, J,$ and $K$ are unit intervals. Leaves of $P$ correspond to sets
  $(a, b, K)$ for fixed $a \in I$, $b \in J$, and $b \neq 0$.
A similar foliation can be defined for truncated prisms.

\bigskip
\begin{figure}[h]
\hspace{.3in} \psfig{figure=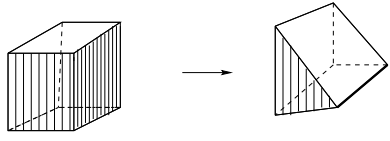, width= 2in} \caption{Induced foliation of a prism} \label{fig2}
\end{figure}
\bigskip

\textbf{\textit{Collapsing}} a (truncated) prism $P$ means taking
the quotient space $P / \sim$, where $x \sim y$ for any two points
, $x$ and $y$, belonging to the same leaf.

We call the \textbf{\textit{sides}} of an $I$-bundle or a tip, the faces which were originally subsets
of the 2-skeleton of $M$. We call the \textbf{\textit{face(s)}} of an $I$-bundle or a tip, the face(s) which
were originally embedded in $S$.

Similarly, we define the \textbf{\textit{sides}}, the \textbf{\textit{top}} face, and the \textbf{\textit{bottom}}
face of a (truncated) prism $P$. Our point here is to give a name for the two hexagonal faces of a truncated prism.
It doesn't matter which is the top and which is the bottom.

\begin{figure}[h]
\hspace{.3in} \psfig{figure=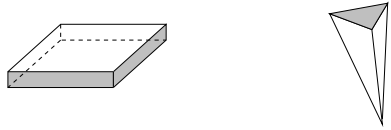, width= 2in} \caption{The sides of an $I$-bundle and the face of a tip}
\label{fig5}
\end{figure}

\begin{figure}[h]
\hspace{.3in} \psfig{figure=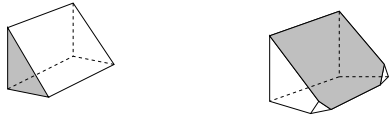, width= 2in} \caption{A side and the top of a (truncated) prism} \label{fig6}
\end{figure}

A polyhedra $P_{1}$ is \textbf{\textit{adjacent}} to a (truncated) prism $P_{2}$ if a side of $P_{2}$ coincides with a
side of $P_{1}$. To each prism, there are at most 2 adjacent polyhedra.

A \textbf{\textit{collapsing annulus}} is an annulus $A$ embedded in the 2-skeleton of $M$ with the following properties:

- $\partial A$ = $\beta_{1}$ $\cup$ $\beta_{2}$ , where $\beta_{1}$ and $\beta_{2}$ are composed of parallel normal arcs.

- $A \cap S = \partial A$.

A \textbf{\textit{collapsing Mobius band}} is a Mobius band $B$
embedded in the 2-skeleton of $M$ with the following
 properties:

- $\partial B$ = $\beta$, where $\beta$ is composed of pairs of
parallel normal arcs.

- $B\cap S = \partial B$.

 A \textbf{\textit{collapsing disk}} is a disk $D$ embedded in the 2-skeleton of $M$ with the following properties:

- $\partial D$ is composed of normal arcs and $D \cap M$ is composed of triangles only.

- $D \cap S = \partial D$.

\bigskip
\begin{figure}[h]
\hspace{.1in} \psfig{figure=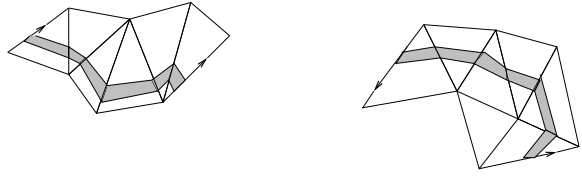, width= 3in} \caption{A collapsing annulus and a collapsing Mobius band in $T^{(2)}.$}
\label{fig8}
\end{figure}
\bigskip

\end{defn}

\begin{rem}
Let $S$ be a non-trivial normal 2-sphere in $M$, and let $A$ be a collapsing annulus. $A$ inherits a
 natural trivial I-bundle structure $S^1 \times I$. Let $D_1$ and $D_2$ be the two disjoint disks on $S$ such that
 $\partial A = \partial D_1 \cup \partial D_2.$  Consider the 3-manifold $M' = M \backslash (S \cup A)$. The
 boundary of $M'$ consist of the union of 2 copies of $A$, $A_1$ and $A_2$,  and 2 (possibly disconnected) copies  of $S$,
 $S_1$ and $S_2$.
 By \textbf{\textit{collapsing}} $A$ in $M'$, we mean taking the quotient space
 $M'$/$(\phi_1, \phi_2, \psi)$ where $\phi_i$ is the natural retraction on $A_i$,
  $\phi_i  : S^{1} \times$ I $\rightarrow $ \{pt\} $\times I$, and where $\psi$ maps any connected
  component of $S_i$ to a point. Note, collapsing $A$ in $M'$ is topologically equivalent to collapsing two
  disjoint
  2-spheres in $M$, parallel to $S$ and $A \cup D_1 \cup D_2$.

\bigskip
\bigskip

A collapsing annulus $A$ is \textbf{\textit{inessential}} if the
2-sphere $D_1 \cup D_2 \cup A$ bounds a ball $B^3$ such that
$int(B^3) \cap S = \emptyset$. Otherwise, $A$ is called essential.
Note, if $A$ is inessential, then collapsing $A$ is topologically
equivalent to collapsing $S$ only.

\smallskip

Let $S$ be a non-trivial normal 2-sphere in $M$, and let $D$ be a
collapsing disk. $D$ inherits a natural foliation as in figure~\ref{fig3.33}
below. Let $D_1$ and $D_2$ be the two disks on $S$, with disjoint
interior, such that
 $\partial D = \partial D_1 = \partial D_2.$ Consider the 3-manifold $M' = M \backslash (S \cup D)$. The
 boundary of $M'$ consists of the union of 2 copies of $D$, $D'$ and $D"$, 1 connected copy of $S$, $S_1$ , and
1 disconnected copy of $S$, $S_2$ . By
\textbf{\textit{collapsing}} $D$ in $M'$, we mean taking the
quotient space
 $M'$/$(\phi_1, \phi_2, \psi)$ where $\phi_1$ (resp. $\phi_2$) is the continuous map which maps each leaf of $D'$
 (resp. $D"$) to a point,
and where $\psi$ maps any connected component of $S_i$ to a point.
Note, collapsing $D$ in $M'$ is topologically
 equivalent to collapsing two disjoint
  2-spheres in $M$ parallel to $S$ and to $D \cup D_1 \cup D_2$.

\bigskip
\begin{figure}[h]
\hspace{.1in} \psfig{figure=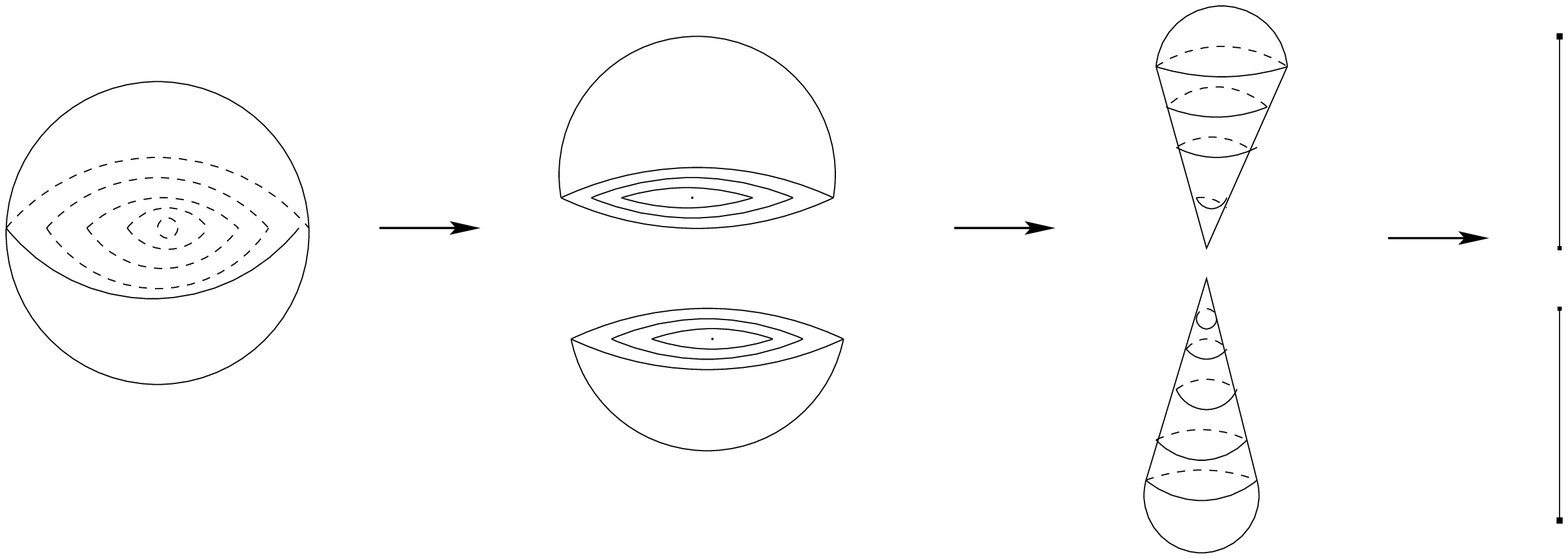, width= 3in} \caption{Collapsing of the disk $D.$}
\label{fig3.33}
\end{figure}
\bigskip

A collapsing disk $D$ is \textbf{\textit{inessential}} if one of the
2-spheres, $D_1 \cup D$ or $D_2 \cup D$,
 bounds a ball $B^3$ such that $int(B^3) \cap S = \emptyset$. Otherwise, $D$ is called essential.
 Note, if $D$ is inessential, then collapsing $D$
  is topologically equivalent to collapsing $S$ only.

\smallskip

Let $S$ be a non-trivial normal 2-sphere in $M$, and let $B$ be a
collapsing Mobius band. The boundary of $B$ is homeomorphic to an
annulus $A$ which inherits a natural trivial $I$-bundle
structure. Let $D_1$ and $D_2$ be the two disks on $S$, with
disjoint interior, such that $\partial D_1 = \partial D_2 =
\partial B$. Let $M' = M \backslash S \cup B$. The boundary of
$M'$ consists of the union of one copy of $A$, one connected copy of $S$, and one disconnected copy of $S$. By
\textbf{\textit{collapsing}} $B$ in $M'$, we mean taking the quotient space $M'$/$(\phi,
\psi)$, where $\psi$ maps any connected component of $S$ to a point and $\phi$ is the natural retraction on $A$.
Collapsing $B$ in $M'$ is topologically equivalent to collapsing two disjoint 2-spheres in $M$, one parallel to $S$
and one parallel to $D_1 \cup D_2 \cup A$ (note, the two 2-spheres bound a twice punctured \R).

\end{rem}

\medskip

A collapsing annulus, Mobius band, or disk is always defined with
respect to an embedded normal surface. These collapsing surfaces will
be a major ingredient in the proof of Theorem~\ref{MT}. We
mention here that a very similar collapsing process has already been
described by Jaco and Rubinstein in ~\cite{JR}.

\section{Proof of Theorem~\ref{MT}}

We are now ready to prove Theorem~\ref{MT}. The proof of this
theorem is presented in the form of a procedure which, given a
normal 2-sphere $S$, finds triangulations of some connected summands
of $M$. The number of resulting summands is completely determined
by $M$ and $S$. Hence, given two topologically parallel
non-isotopic normal 2-spheres, we could obtain different
decompositions of $M$.


Let us start by cutting $M$ along $S$. We obtain a cell
decomposition of $M \backslash S$ composed of the 7 types of
polyhedra described earlier. We summarize the procedure: Let
$S$ be the non-trivial normal 2-sphere which we collapse. We
get rid of all tips and $I$-bundles by collapsing some collapsing annuli
and disks. Each maximal collection of adjacent tips contributes to
an $S^3$ summand and each maximal collection of adjacent
$I$-bundles contributes to either an $S^3$ or an \R summand. We
then collapse the prisms one at a time. Each collapsing may or may
not contribute to one $L(3, 1)$, one or two \R, or one $S^3$
summand. Finally, we collapse each truncated tetrahedra and we
count the number of $S^1 \times S^2$.

{\bf Step 1:} Consider a maximal collection
$\mathcal{C}$ of adjacent tips. First, note that a tip can only be adjacent to
another tip or a prism.  The union of the faces of these
tips represents a connected subsurface $S'$ of $S$, and so
$\mathcal{C}$ can be seen as the cone over $S'$. Call $\partial
\mathcal{C}$ the union of the sides along which a tip in $\mathcal
{C}$ is adjacent to a prism. By definition, $\partial
\mathcal{C}$ is is made of a nonempty union of collapsing disks. We collapse each such disk. Note,
the collapsing of each collapsing disk preserves the cone structure
of $\mathcal{C}$. In particular, each connected boundary component
of $S'$ is being collapsed to a point, and so $S'$ is being
collapsed to a 2-sphere. Hence, for each maximal collection
$\mathcal{C}$ of adjacent tips, we obtain a cone over a 2-sphere
which is itself homeomorphic to a ball. It is now clear that each
$\mathcal{C}$ contributes to a trivial summand in the
decomposition of $M$.

\bigskip
\begin{figure}[h]
\hspace{.3in} \psfig{figure=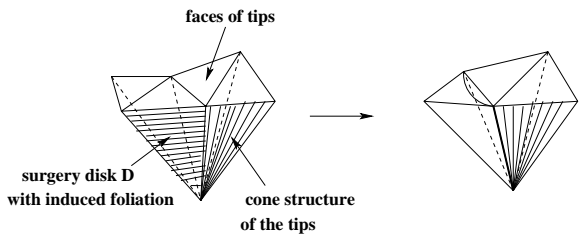, width= 3in} \caption{The collapsing of a collapsing disk along the side of a
 tip preserves the cone structure of $\mathcal{C}.$} \label{fig3.2}
\end{figure}
\bigskip

{\bf Step 2:} Consider a maximal collection
$\mathcal{J}$ of adjacent $I$-bundles. Note, an $I$-bundle can only be adjacent to another
$I$-bundle or to a truncated prism. The union of the faces of
these $I$-bundles represents a possibly disconnected subsurface
$S'$ of $S$, and so $\mathcal{J}$ can be seen as an $I$-bundle
over a surface $B$ whose double cover is $S'$. Call $\partial
\mathcal{J}$ the union of the sides along which an $I$-bundle in
$\mathcal{J}$ is adjacent to a truncated prism. If $\partial
\mathcal{J}$ is empty then $\mathcal{J}$ is an $I$-bundle whose
boundary is $S$ which means that $\mathcal{J}$ is homeomorphic to
a punctured \R. If $\partial \mathcal{J}$ is nonempty, then it
consists of a union of collapsing annuli. We
collapse each such annulus. Note, the collapsing of each
annulus preserves the $I$-bundle structure of $\mathcal{J}$. In
particular, each connected boundary component of $S'$ is being
collapsed to a point, and so $S'$ is being collapsed to either one
or two copies of a 2-sphere depending on whether $S'$ is
disconnected or not. Hence, for each maximal collection
$\mathcal{J}$ of adjacent prisms, we obtain an $I$-bundle whose
boundary is either one or two copies of a 2-sphere. It is now
clear that each $\mathcal{J}$ contributes to either a trivial
summand or an \R \space summand in the decomposition of $M$.

\bigskip
\begin{figure}[h]
\hspace{.3in} \psfig{figure=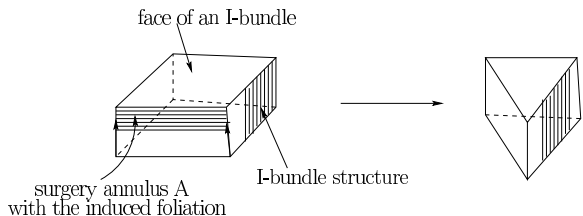, width= 3in} \caption{The collapsing of a collapsing annulus along the side of an
$I$-bundle preserves the $I$-bundle structure of $\mathcal{J}.$} \label{fig3.1}
\end{figure}
\bigskip

\begin{rem}
Let $P$ be a (truncated) prism adjacent to a tip or an $I$-bundle $Q$. From step 1 and 2, we know that
the common side $F$ of $P$ and $Q$ is part of a collapsing annulus or disk $A$. Moreover, the collapsing of $A$
induces a collapsing of $F$. Note, this collapsing of $F$ coincides with our definition of the collapsing of $P$.
Hence, after collapsing all the collapsing surfaces in step 1 and 2, we obtain a cell decomposition  of $M \backslash
S$ made of (truncated) tetrahedra, (truncated) prisms, and (truncated) prisms with their sides collapsed
(these cells are called pillows in ~\cite{JR}).
By abuse of language, we also call the latter cells prisms.

\end{rem}

{\bf Step 3:} This step consists of collapsing each (truncated) prism one at a time.
Let $P$ be a (truncated) prism. If $P$ is embedded in $M \backslash S$, we can
collapse it without changing the topology of $M \backslash S$. If $P$ is not embedded, there are 4 cases to
consider.

\noindent {\bf Case 1:} All the leaves of one (or both) side of
$P$ are not embedded in $M$, and all the other leaves of $P$ are.
This happens if one (or both) side of $P$ represents a collapsing
annulus or disk depending on whether $P$ is truncated or not. We
cut along the collapsing surface and we collapse $P$.

Each side of $P$ may be embedded but their union may represent a
collapsing annulus or disk. In this case we cut along both sides and we collapse.

\bigskip
\begin{figure}[h]
\hspace{.3in} \psfig{figure=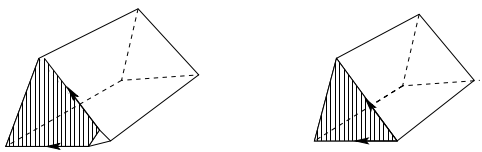, width= 2.5in} \caption{A collapsing annulus and a collapsing disk.} \label{fig19}
\end{figure}
\bigskip

\noindent {\bf Case 2:} Exactly one leaf in one (or both) side of
$P$ is not embedded, and all the other leaves are. This happens when
 one (or both) side of $P$ represents a collapsing Mobius band
$B$. We cut along $B$, we collapse $P$, and we add one (or two) \R summand in the
decomposition of $M$.

Each side of $P$ may be embedded but their union may represent a
collapsing Mobius band $B$. In this case, we cut along both sides, we collapse $P$,
 and we obtain one \R summand in the decomposition of $M$.

 \noindent {\bf Case 3:} None of the leaves of $P$ are embedded.
 This happens when the top of $P$ is identified to its bottom
 without any twist. In particular, both sides of $P$ are collapsing
 annuli or disks depending on whether $P$ is truncated or not. We
 collapse the collapsing surfaces. One can check that after the collapsing, $P$ is homeomorphic to $S^3$.
So we simply remove $P$ from the cell decomposition of $M \backslash S$.

\bigskip
\begin{figure}[h]
\hspace{.3in} \psfig{figure=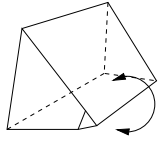, width= .8in} \caption{A prism with no embedded leaves.} \label{fig3.10}
\end{figure}
\bigskip

\noindent {\bf Case 4:} Exactly one leaf of $P$, which does not
belong to one of its sides, is not embedded. This happens when the
top and bottom of $P$ are identified by a 1/3 twist. Here, $P$ is
homeomorphic to a solid torus $T$ whose boundary consists of a
collapsing annulus $A$ and an annulus $A' \subset S$. Since $A'$ lies
on $S$ both of its boundary components bound disjoint disks $D_1$
and $D_2$ on $S$. Moreover, both $\partial D_1$ and $\partial D_2$
represent a (3,1) curve (or (3, 2) curve) on $\partial T$. Hence
$Nbhd(T \cup D_1 \cup D_2)$ is homeomorphic to a twice punctured lens
space $L(3, 1)$. We remove $P$ from the cell decomposition of $M \backslash$. We then collapse $A$ and we add one
 $L(3, 1)$ summand in the decomposition of $M$.

\bigskip
\begin{figure}[h]
\hspace{.3in} \psfig{figure=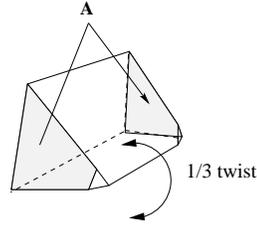, width= 1.3in} \caption{A truncated prism whose boundary is homeomorphic to
a torus.} \label{fig3.13}
\end{figure}
\bigskip

\begin{figure}[h]
\hspace{.3in} \psfig{figure=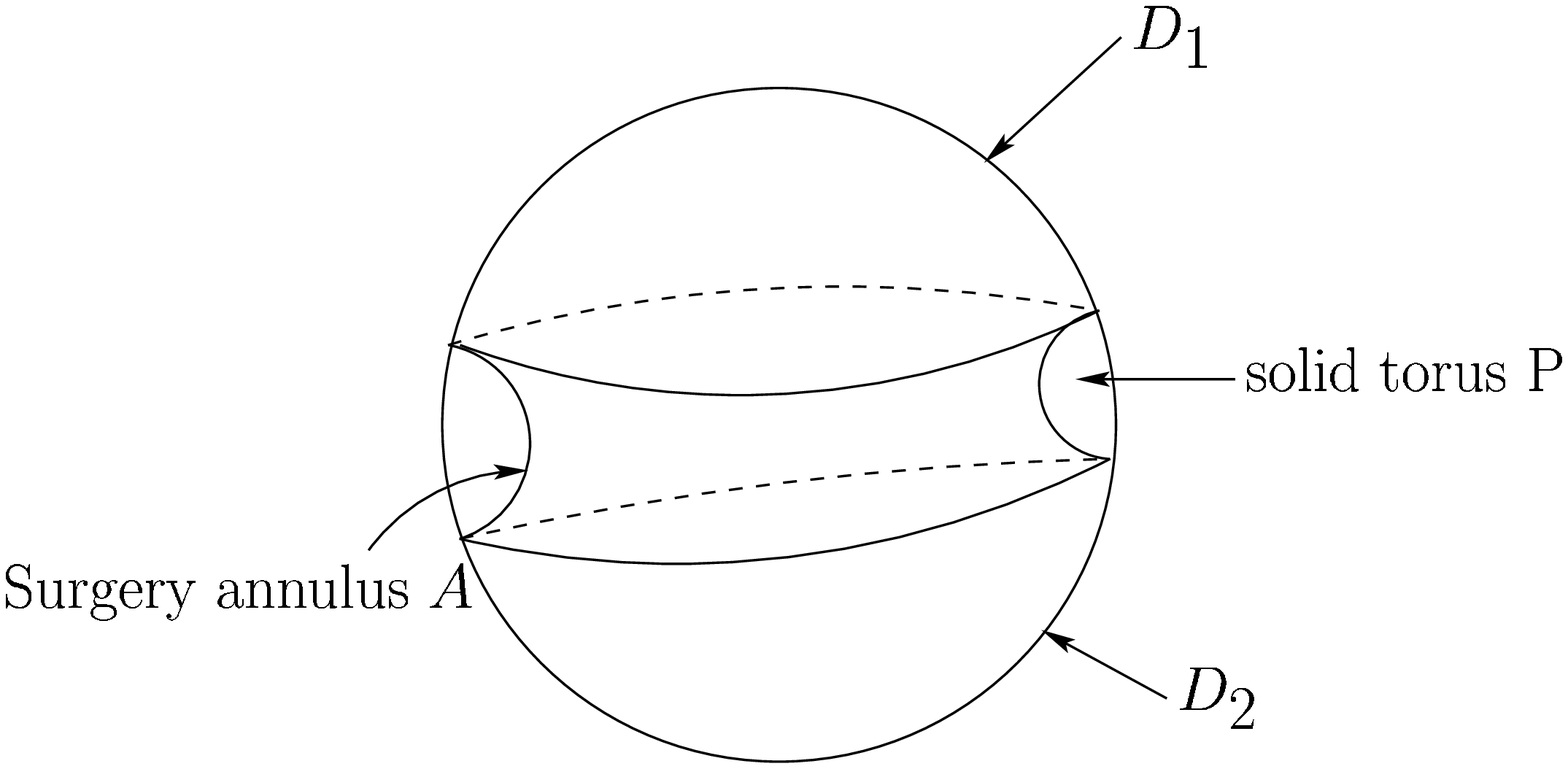, width= 2in} \caption{A twice punctured $L(3, 1)$.} \label{fig3.17}
\end{figure}
\bigskip

Suppose that the prism $P$ was adjacent to another prism $P'$.
Then the collapsing of $P$ induced a collapsing of one (or both) side of $P'$.
Since the collapsing of two adjacent prisms coincide along their adjacent sides, $P'$ is thought as a prism with
 one (or both) of its sides collapsed (these cells are called pillows in ~\cite{JR}). Treating  each $P'$ as
such, we repeat step 3 for each prism (with or without their sides collapsed) found in the new
cell decomposition of $M$.

{\bf Step 4:} The only polyhedra left are tetrahedra and truncated tetrahedra. We collapse each truncated
tetrahedron to a regular tetrahedron.

{\bf Step 5:} We need to count the number of $S_1 \times S^2$
summands. Recall that collapsing a collapsing disk or annulus is
topologically equivalent to collapsing a 2-sphere. Let $n$, $m$,
and $p$ be the respective number of collapsing annuli, disks and
Möbius bands found in step 1, 2, and 3. Let $k$ be the number of
connected summands (including the trivial ones) of $M$ obtained in
step 1, 2, and 3. Then the number of $S^1 \times S^2$ summands is
$n + m + p + 2 - k$. This comes from the fact that each
non-separating 2-sphere contributes to one such summand. Indeed,
if we have collapsed $n + m + p$ collapsing surfaces, we have
actually collapsed $n + m + p + 1$ 2-spheres.

Note, if $M_1$, ..., $M_k$ are the resulting triangulated
summands of $M$, then it follows directly from the above
construction that $|M_1| +$ ... $+ |M_k| < |M|$. In particular, if
$S$ has k non-zero quadrilateral types, then  $|M_1| +$ ... $+
|M_k| = |M| - k.$

\section{Connected sums of triangulated 3-manifolds}
\label{consum}

We describe here geometrical constructions to obtain a
triangulation for the connected sums of closed orientable
triangulated 3-manifolds. Intuitively, if one wants to take the
connected sum of two such 3-manifolds, one can think of removing the
interior of a tetrahedron in each manifold and glue the resulting
manifolds along their boundaries. The problem is that
the boundary of a tetrahedron may not be homeomorphic to a
2-sphere. One could always retriangulate a tetrahedron in each
manifold to obtain tetrahedra with embedded boundary, but this construction seems to be artificial and
 is not as efficient as the one described here. Construction~\ref{cons1} is the only construction which is
 described thoroughly. As the reader will become more familiar with it, the other constructions are simple
 generalizations of the first one.

\smallskip

\begin{construction} \label{cons1}
Let $P$ and $N$ be two triangulated closed orientable 3-manifolds
with $t_{1}$ and $t_{2}$ tetrahedra, and $v_{1}$ and $v_{2}$
vertices, respectively. If not both of $v_{1}$ and
$v_{2}$ are equal to 1, then there is a triangulation of $P \# N$
with $t_{1} + t_{2} + 2$ tetrahedra and $v_{1} + v_{2} - 2$
vertices. If $v_{1} = v_{2} = 1$, then there is a 1-vertex
triangulation of $P \# N$ and $t_{1} + t_{2} + 4$ tetrahedra.
\end{construction}

\medskip

This triangulation of $P \# N$ does not have to be minimal, even
if $P$ and $N$ are minimal. There are some cases, depending on the
triangulations of P and N, where only $t_{1} + t_{2} + 1$ tetrahedra are needed to construct $P \# N$.

We assume that either $P$ or $N$ has more
than one vertex. Indeed, let $M$
be the one-vertex triangulation of a closed orientable 3-manifold
which is the connected sum of exactly two irreducible 3-manifolds
$P$ and $N$, and let $S$ be an essential normal 2-sphere.
We collapse $S$. It is shown in Lemma~\ref{L4} that one summand,
 say $P$, has one vertex (from the collapsing of $S$) and that $N$ has two vertices (one from the collapsing of $S$ and
  one from the original triangulation of $M$). Hence, it seems somewhat natural to make this assumption.

We first assume that $v_{1}$ and $v_{2}$ are both strictly greater than 1. Let $a$ be a vertex of $P$ and $b$ a
vertex of $N$. We remove a normal neighborhood of $a$ and $b$. $P \backslash S_{1}$ and $N \backslash S_{2}$ are
 now 3-manifolds composed of tetrahedra and truncated tetrahedra with a boundary component being a triangulated
 2-sphere, $S_{1}$ and $S_{2}$ respectively. We want to change the cell decompositions of $P \backslash S_{1}$
  and $N \backslash S_{2}$ in order to glue the two manifolds along their boundary and obtain a well-defined
  triangulation for their connected sums. To simplify the notation, we will denote $P \backslash S_{1}$ and $N
  \backslash S_{2}$ by $P'$ and $N'$, respectively.

Because $P$ and $N$ have more than one vertex, we assume without loss of generality that $a$ and $b$
are chosen so that $P'$ and $N'$ contain a face of a truncated tetrahedron as in figure~\ref{fig4.1} below.

\bigskip
\begin{figure} [h]
\centering{\psfig{figure=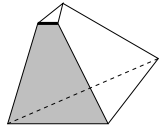, width= .8in}

   \caption{A
truncated tetrahedron in $P'$ and $N'$.} \label{fig4.1}}
\end{figure}
\bigskip

For future reference, we keep track of the above thickened edges
in the triangulations of $S_{1}$ and $S_{2}$.

\begin{defn}

We say that a face of a tetrahedron in a triangulation is a
\textbf{\textit{cone}} if two of its edges are identified as in
the figure~\ref{fig4.2} below. If, in addition to that, the endpoints of each edge are distinct, we call the cone a
\textbf{\textit{good cone}}.

\bigskip
\begin{figure}[h]
\centering{\psfig{figure=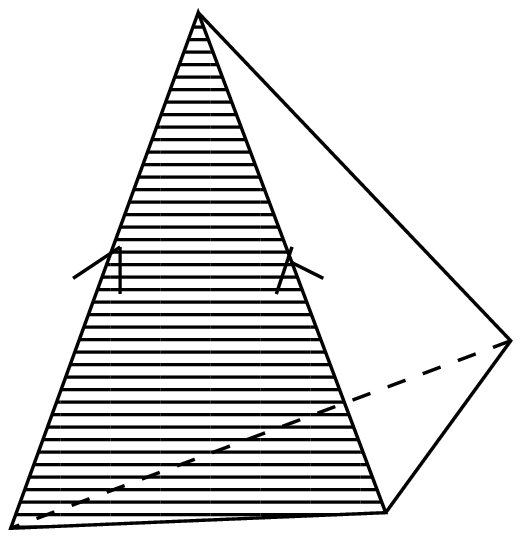, width= 0.7in} \caption{A cone
in a tetrahedron.} \label{fig4.2}}
\end{figure}
\bigskip

Let $\tau$ be a triangulation of a surface $S$. The \textbf{\textit{dual}} $\tau'$ of $\tau$ is a cell
decomposition of $S$ such that :

-There is a 1-1 correspondence between the i-simplices of $\tau$ and the $(2 - i)$-simplices of $\tau'$.

-Every triangle of $\tau$ (resp. every polygon of $\tau'$) contains exactly one vertex of $\tau'$ (resp. $\tau$).

-Every edge of $\tau$ (reps. $\tau'$) intersects exactly one edge of $\tau'$ (resp. $\tau$) in exactly one point.

\end{defn}

One can check that for any fixed triangulation there is
exactly one dual and vice versa. We denote by \textbf{\textit{G}}
the 3-valent planar graph made of the union of the vertices and
edges of $\tau'$. Let $T_{1}$ and $T_{2}$ be the triangulations
corresponding to the 2-spheres $S_{1}$ and $S_{2}$, respectively.
Let $G_{1}$ and $G_{2}$ be the respective graphs corresponding to
the respective duals $\tau'_1$ and $\tau'_2$ of $T_{1}$ and
$T_{2}$.

Suppose we have colored the vertices of $\tau'_1$ and $\tau'_2$.
We say that $G_1$ and $G_2$ are \textbf{\textit{colored
homeomorphic}} (we write $G_1 = \overline{G_2}$) if there exists a
cell preserving orientation reversing homeomorphism of the
2-sphere, sending $\tau'_1$ to $\tau'_2$ and sending each colored
vertex of $\tau'_1$ to a same color vertex in $\tau'_2$. We
describe a set of rules to transform the graphs $G_{1}$ and
$G_{2}$ into graphs $G'_1$ and $G'_2$ such that $G'_1$ is colored
homeomorphic to $G'_2$.

Let $G_{1}$ and $G_{2}$ be given.
We first color the vertices of $G_{1}$ in white and the vertices of $G_{2}$ in black.
 We can add 4-valent red vertices on the original edges,
 and 3-valent black and white vertices and edges using the following rules:

\begin{enumerate}
\item - Black vertices cannot be joined by an edge to white vertices.

\item - Black (resp. white) vertices cannot be added to $G_{2}$ (resp. $G_{1}$). (This rule will be omitted when
either $P$ or $N$ has a one-vertex triangulation).

\item - Let $r$ be a red vertex. Let $e_{1}$ and $e_{2}$ be
opposite edges with one common end point $r$ (since r is 4-valent,
the notion of opposite edges is well-defined). Let $a_{1}$ and
$a_{2}$ be the two other endpoints of $e_{1}$ and $e_{2}$,
respectively. If $a_{1}$ is black (resp. white), then $a_{2}$
cannot be white (resp. black).

 \item - The
resulting graphs $G'_{1}$ and $G'_{2}$ must be connected and
planar.
\end{enumerate}

\bigskip
\begin{figure}[h]
\centering{\psfig{figure=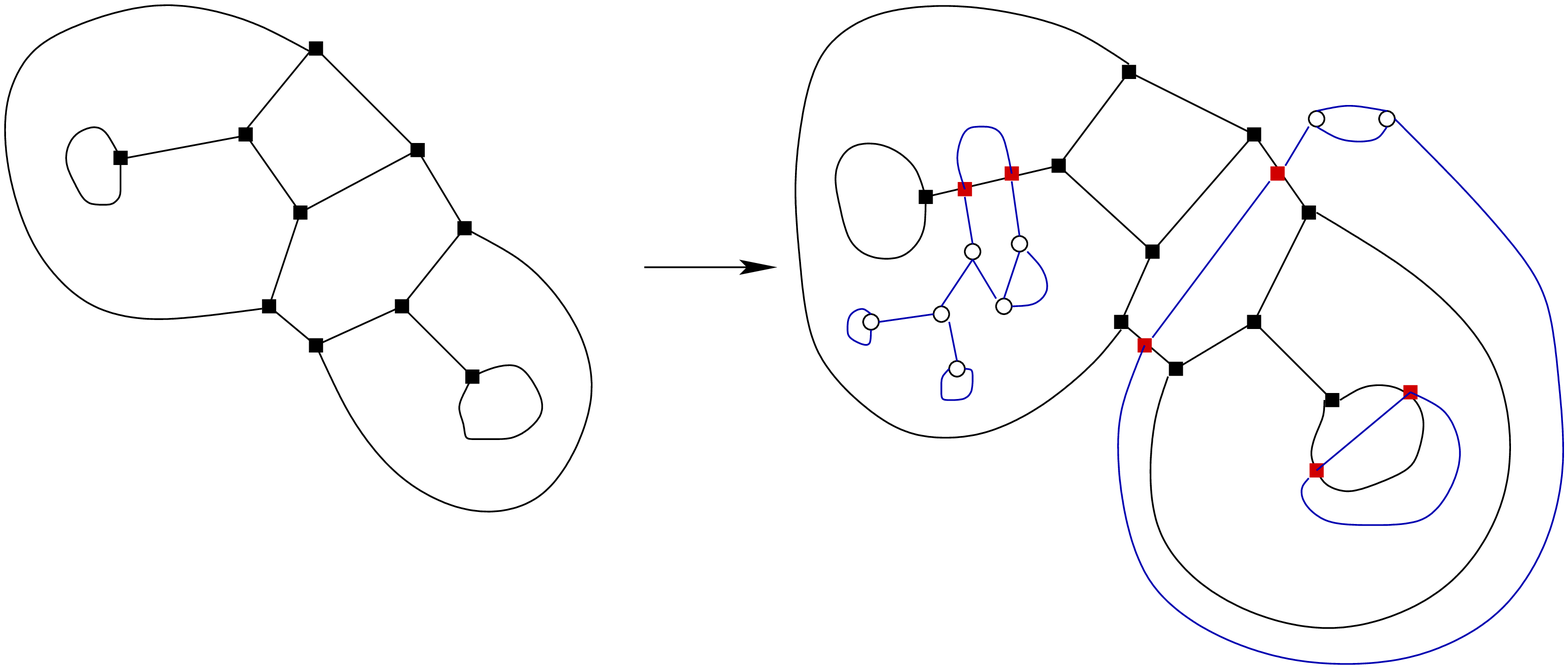, width= 3in} \caption{Example of a 3-valent graph $G$ and another graph $G'$
obtained by adding edges and vertices according to the above rules.} \label{fig4.4}}
\end{figure}
\vspace{.2in}

Let $S_{1}$ and $S_{2}$ be given. To the thickened edge in
figure~\ref{fig4.1} corresponds a unique edge in the induced
graphs $G_{1}$ and $G_{2}$, say $e_{1}$ and $e_{2}$ respectively.
We insert two red vertices on each of these edges: $r_{1}$ and
$r'_{1}$ on $G_{1}$, and $r_{2}$ and $r'_{2}$ on $G_{2}$. We then
connect $r_{i}$ and $r'_{i}$ by an edge. Consider now a copy of
the colored homeomorphic image of $G_{2}$ with the edge $e_{2}$
removed. We call this new graph \GGG\ . Color all the vertices of
\GGG\ in black. Draw an edge emanating from $r_{1}$ and one
emanating from $r'_{1}$. On these two edges, draw the graph
corresponding to \GGG\ . We call the resulting graph $G'_1$. The same procedure can be done by adding
\GG\ to $G_{2}$ to obtain the graph $G'_2$. One can check that $G'_{1}$ is colored
homeomorphic to $G'_2$.

\bigskip
\begin{figure}[h]
\centering{\psfig{figure=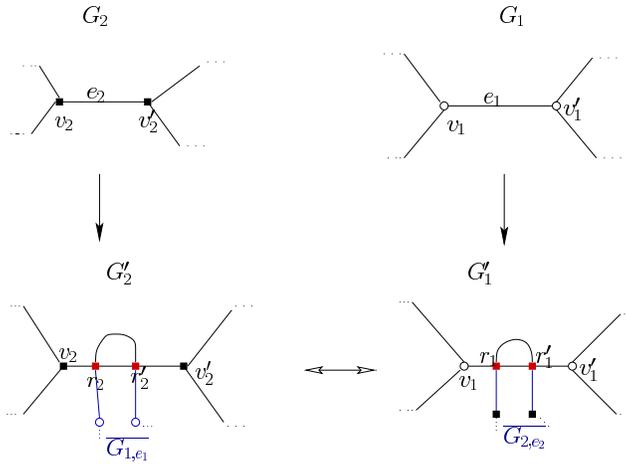, width= 3.3in}
\caption{Colored homeomorphism between $G'_1$ and $G'_2$.}
\label{fig4.5}}
\end{figure}
\vspace{.2in}

Let us take a closer look at the triangulations of $S^1$ and $S^2$. Note, because $G'_{1}$ and $G'_{2}$ are planar
and connected, they represent the cell decompositions of 2-spheres.
  This implies that their duals, $S'_{1}$ and $S'_{2}$, also represent the cell decompositions of 2-spheres.
  From the definition of a dual triangulation, each 3-valent (resp. 4-valent) vertex added on $G_{1}$ corresponds to
   a triangle (resp. quadrilateral) added on $S_{1}$, and each edge added on $G_{1}$ corresponds
   to an edge added on $S_{1}$.
  Hence, we changed the triangulation of $S_{1}$ and $S_2$ into cell decompositions of 2-spheres with
  triangles and quadrilaterals.

\bigskip
\begin{figure}[h]
\centering{\psfig{figure=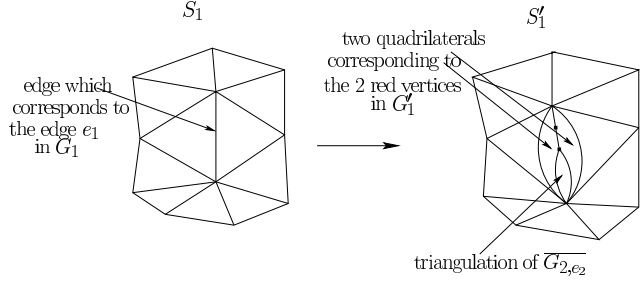, width= 3.3in} \caption{View of the new cell decomposition of $\partial P'$.}
\label{fig4.6}}
\end{figure}
\vspace{.2in}

 We saw that each 4-valent vertices added on $G_{1}$ corresponds to adding a quadrilateral on $S_{1}$.
 Here, adding a quadrilateral on $S_{1}$ is represented by inserting a prism in $P'$, and adding a
 triangle on $S_{1}$ is represented by inserting the tip of a tetrahedron in
 $P'$. See figure~\ref{fig4.7} below.

As described above, the black vertices that we added on $G_{1}$, which correspond to tips of tetrahedra inserted in $P'$,
 represent a graph isomorphic to $\overline{G_{2, e_{2}}}$.
 Since $G_{2}$ is the dual of the triangulation of a 2-sphere, $\overline{G_{2, e_{2}}}$ is the dual of the
  triangulation of a disk.
  Hence, the union of all the tips inserted in $P'$ is homeomorphic to the cone over a disk, which is homeomorphic
  to a 3-ball.
  Hence, inserting the two prisms and the tips of tetrahedra did
   not change the topology of $P'$ and $N'$. Moreover, every white (resp. black) triangle in $S'_{1}$
  (resp. $S'_2$) corresponds to a
  truncated tetrahedron, and every white (resp. black) triangle in $S'_{2}$ (resp.
  $S'_1$)
  corresponds to the tip of a tetrahedra. If $\phi$ is the colored
  homeomorphism from $G'_1$ to $G'_2$, it is now clear that $(P \cup N)$/$\phi$
   defines a well-defined triangulation for $P\# N$.

\begin{figure}[h] ~\label{fig4.7}
\centering{\psfig{figure=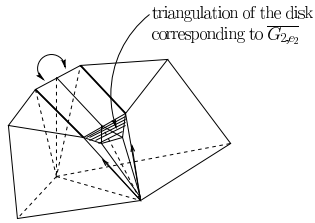, width= 1.6in} \caption{The new cell decomposition of $P'$.} \label{fig4.7}}
\end{figure}

\bigskip
\begin{figure}[h]
\centering{\psfig{figure=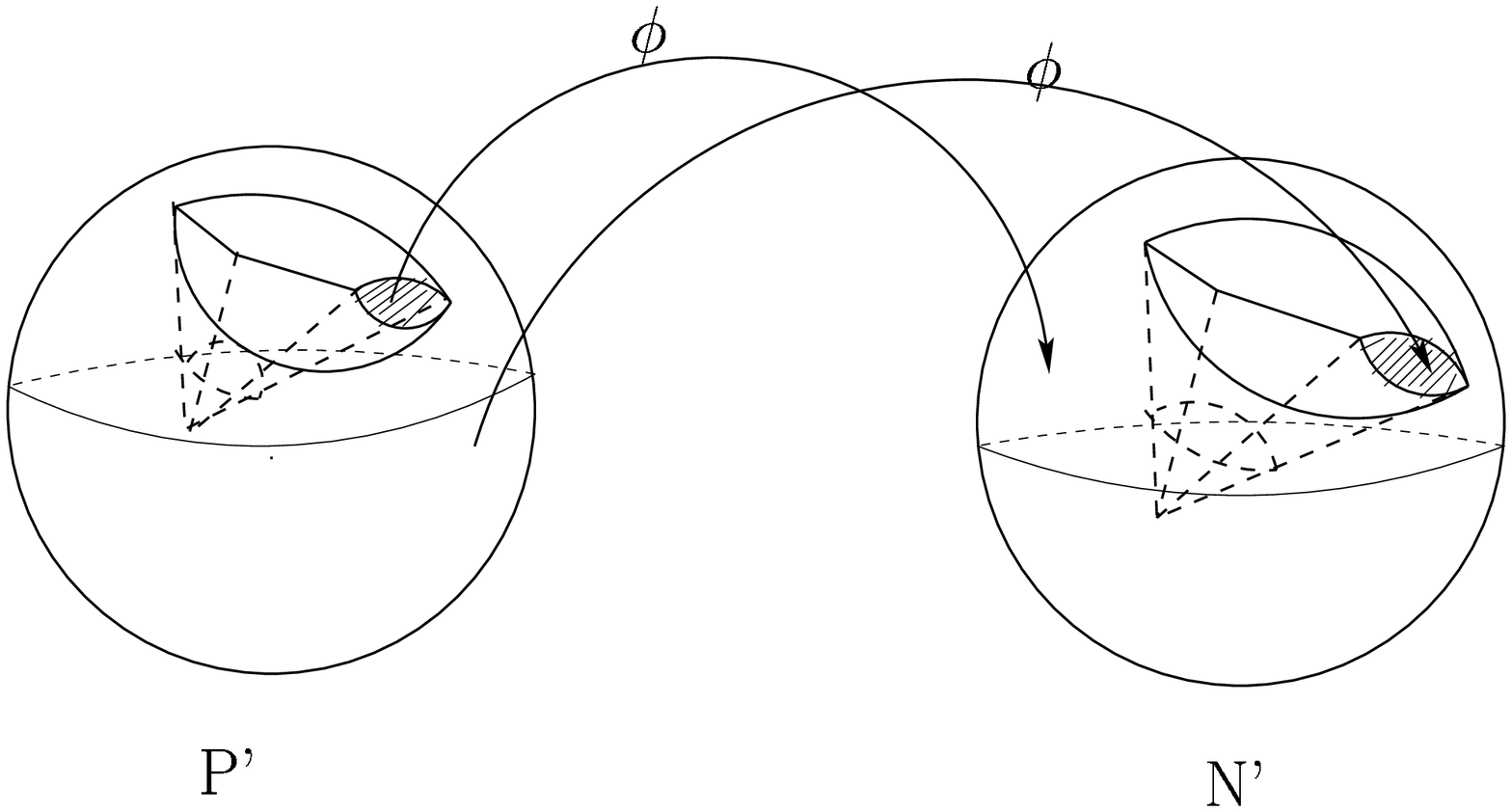, width= 2.4in} \caption{The connected sum of $P$ and $N$.} \label{fig4.8}}
\end{figure}
\vspace{.2in}

\hspace{.3cm} Suppose now that $N$ has more than one vertex and
that $P$ has exactly one vertex. Let $G_{1}$ and $G_{2}$ be the
induced graphs of $P$ and $N$ respectively. Consider $G_{1}$. As
before, we add two 4-valent vertices which correspond to inserting
two truncated prisms. Because the prisms are truncated, each time
we add a 4-valent vertex on $G_{1}$ we need to add a pair of white
vertices too. Note, we are not given a choice on where to
place this pair of white vertices. Note also that each black
vertex added on $G_1$ corresponds to a truncated tip. Hence, when
we add black vertices on $G_{1}$ corresponding to the graph
$\overline{G_{2, e_{2}}}$, we also add white vertices on $G_1$
which correspond to the graph of $G_{2, e_{2}}$.
 After adding all the necessary vertices on $G_{1}$, we obtain a graph as in figure~\ref{fig4.9} below.

\bigskip
\begin{figure}[h]
\centering{\psfig{figure=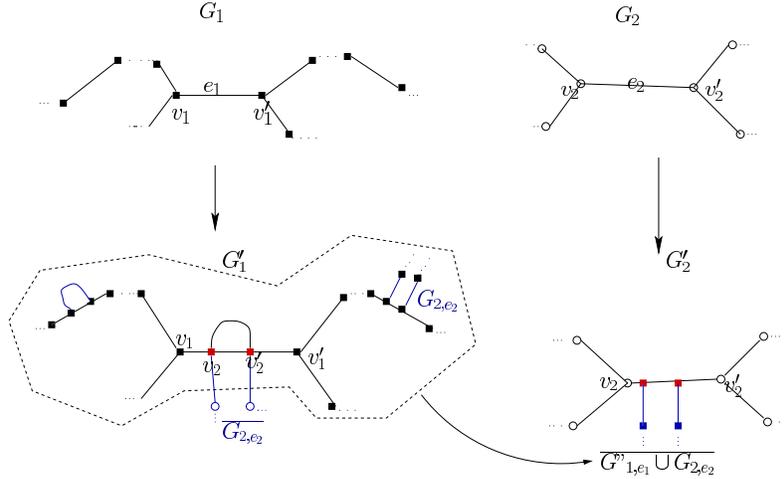, width= 4.1in}
\caption{Colored homeomorphism between $G'_1$ and $G'_2$.}
\label{fig4.9}}
\end{figure}
\vspace{.2in}

 Suppose now that both $P$ and $N$ have exactly one vertex in their triangulations.
 If we now remove a normal neighborhood of each vertex, we end up with a cell decomposition with no vertices
  which is impossible.
  On the other hand, we can take a subdivision of a tetrahedron of, say, $P$, to obtain a
  new triangulation of $P$ with two vertices (see figure~\ref{fig4.10} below).
  We then apply the construction above to the new triangulation of $P$ and to $N$. The triangulation
  of $P \# N$ has $|P_{new}| + |N| + 2 = (|P| + 3) + |N| + 2 = |P| + |N| + 5$ tetrahedra.

\bigskip
\begin{figure}[h] \label{fig4.10}
\centering{\psfig{figure=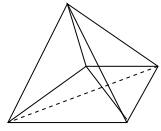, width= .8in} \caption{A
subdivision of a tetrahedron in $P$.} \label{fig4.10}}
\end{figure}
\vspace{.2in}

Let us describe a construction which only requires four more tetrahedra. Consider the 2-tetrahedron 3-vertex
triangulation of $S^{3}$ in figure~\ref{fig4.11} below. Take the connected sum of $P$ with $S^{3}$ by removing
normal neighborhoods of the vertex $v_1$ in $S^{3}$, and $v$ in $P$.
 This vertex $v_1$ as the nice property that its induced 3-valent graph $G$ is 1-edge-connected.
 For future reference, we will say that $v_1$ is a \textbf{\textit{good vertex}}.
 It is not hard to see that a vertex $A$ is a good vertex if it is the end point of two edges forming a good cone.

\bigskip
\begin{figure}[h]
\centering{\psfig{figure=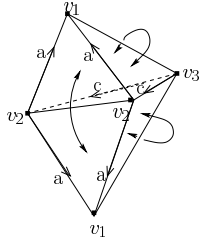, width= 1in}  \caption{A
3-vertex 2-tetrahedron triangulation of $S^3$.} \label{fig4.11}}
\end{figure}
\vspace{.2in}

$P \# S^{3}$ is made of the tetrahedra from $P$, the tetrahedra from $S^{3}$, plus one extra tetrahedron. In fact,
when there is a good vertex in the triangulation of one of two 3-manifolds, there is only one red vertex
 that needs to be added to $G_1$ and $G_2$ to obtain an orientation reversing isomorphism between $G'_1$ and $G'_2$.
Hence, we obtain a triangulation of $P$ with 3 more tetrahedra and 2 vertices.
What we need to notice here is that one of the two vertices of the new triangulation of $P$ is also a good vertex
 (vertex $v_3$ in figure~\ref{fig4.11}).
Hence, $|P_{new} \# N| = |P_{new}| + |N| + 1 = |P| + 3 + |N| + 1 = |P| + |N| + 4$.

\bigskip
\begin{figure}[h]
\centering{\psfig{figure=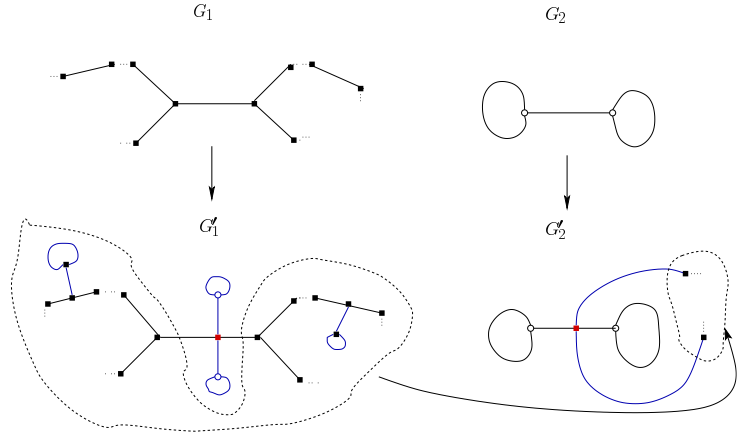, width= 3.8in} \caption{Only one red vertex is needed to obtain $G'_1$ and $G'_2$.}
\label{fig4.12}}
\end{figure}
\vspace{.2in}

 We summarize the above constructions.
 If $P$ has two vertices, then we can take the connected sum of $P$ and $N$ by inserting 2  tetrahedra unless $P$
 has a good vertex in which case the insertion of only 1 tetrahedron is necessary.
 If both $P$ and $N$ have one vertex, then $P \# N$ requires inserting 4 tetrahedra (we will see later that if either
 $P$ or $N$ has a good vertex then this construction requires the insertion of 2 tetrahedra only).

\begin{construction} \label{cons2}
Let $M_{1}$, ..., $M_{k}$ be closed orientable 3-manifolds with at least 2 vertices in each
 of their triangulations.
Then there exists a triangulation for $M_{1} \#$ ...$\# M_{k}$ with $\sum |M_{i}| + k + 2$
 tetrahedra and at least 2 vertices.
\end{construction}

\begin{construction} \label{cons3}
Let $M$ be a triangulated closed orientable 3-manifold with $|M| =
n$ and with exactly one vertex. Then, for $k \geq 3$, there exists
a $k$-vertex triangulation for $M$ with $(n + k + 1)$ tetrahedra.
\end{construction}

\begin{construction} \label{cons4}
Let $M$ be a triangulated closed orientable 3-manifold with $|M| =
n$ and with exactly $k$ vertices. Then, there exists a
$(k+1)$-vertex triangulation of $M$ with $(n + 2)$ tetrahedra. If
$M$ has a cone, this construction can be done by adding a single
tetrahedron, and hence, obtaining a $(k+1)$-vertex triangulation
of $M$ with (n + 1) tetrahedra.
\end{construction}

\begin{construction} \label{cons5}
Let $M_{1}$, ..., $M_{k}$ be closed orientable 3-manifolds with 1-vertex triangulations, and $k \geq 3$. Then there
exists a triangulation for $M_{1} \#$ ...$\# M_{k}$ with $\sum |M_{i}| + 2k$ tetrahedra and 1 vertex.
\end{construction}

We show in ~\cite{Bar2} that construction~\ref{cons5} is the most
efficient way of taking connected sums.

\begin{construction} \label{cons6}
Let M be a triangulated closed orientable 3-manifold with $|M| =
n$ and at least two vertices. Then there exists a triangulation
for $M \# L(3, 1)$ with (n + 2) tetrahedra.
\end{construction}

Let $M'$ be the manifold obtained from $M$ after removing the link
of a vertex. Without loss of generality, we assume $M'$ contains a
truncated tetrahedra as in figure~\ref{fig4.1}. Let $T$ be the
triangulation of $\partial M'$, and $G$ its dual. Let $e$ be the
dual of the thickened edge from figure~\ref{fig4.1}. We insert 2
prisms along the shaded face, and we insert tips of tetrahedra so
that the new dual $G'$ of $\partial M'$
 correspond to figure~\ref{fig4.28} below.

\begin{figure}[h]
\centering{\psfig{figure=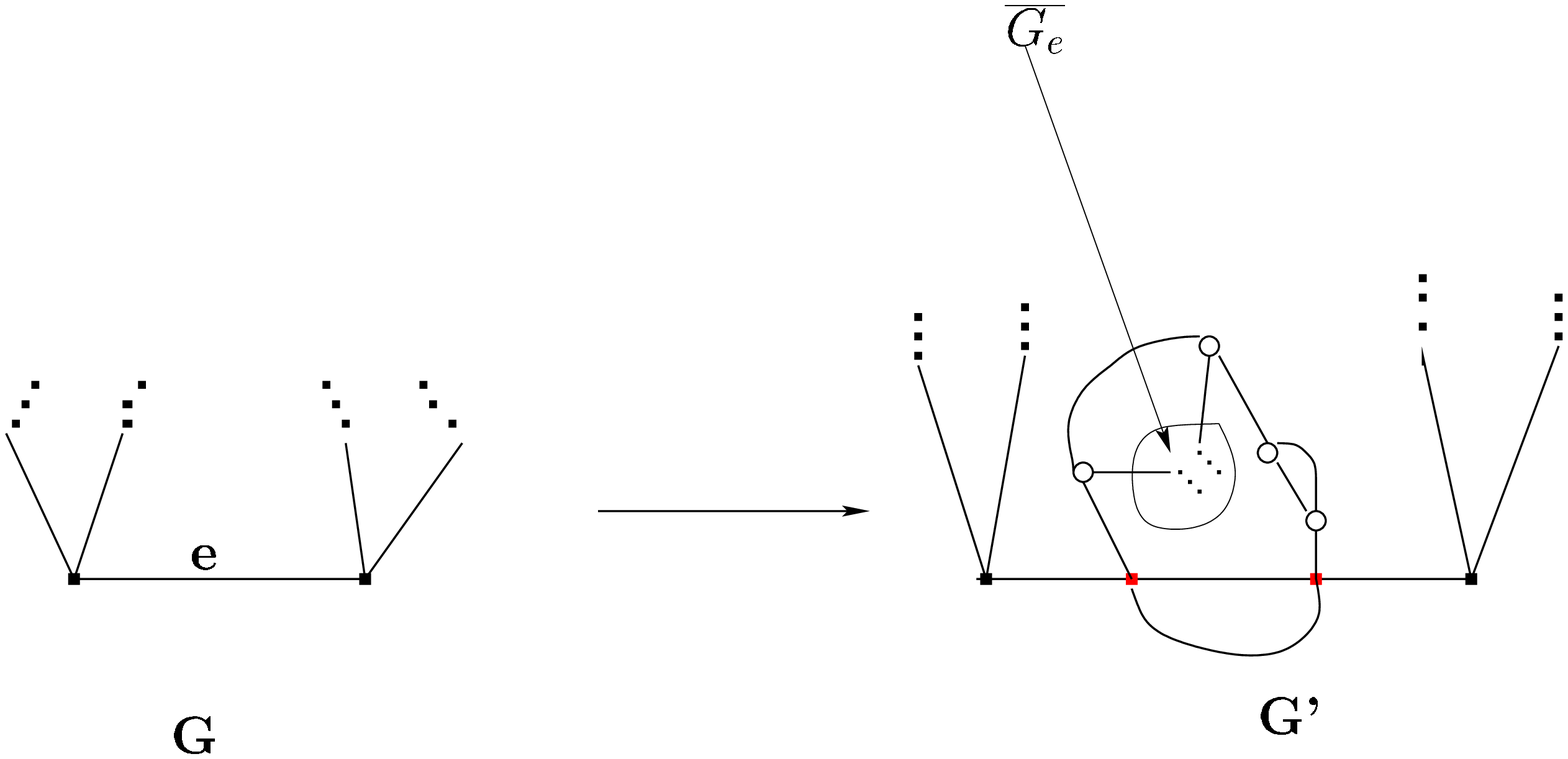, , height= 1in} \caption{The
dual, $G'$, of the new cell decomposition of $\partial M'$.}\label{fig4.28}}
\end{figure}

\bigskip
\begin{figure}[h]
\centering{\psfig{figure=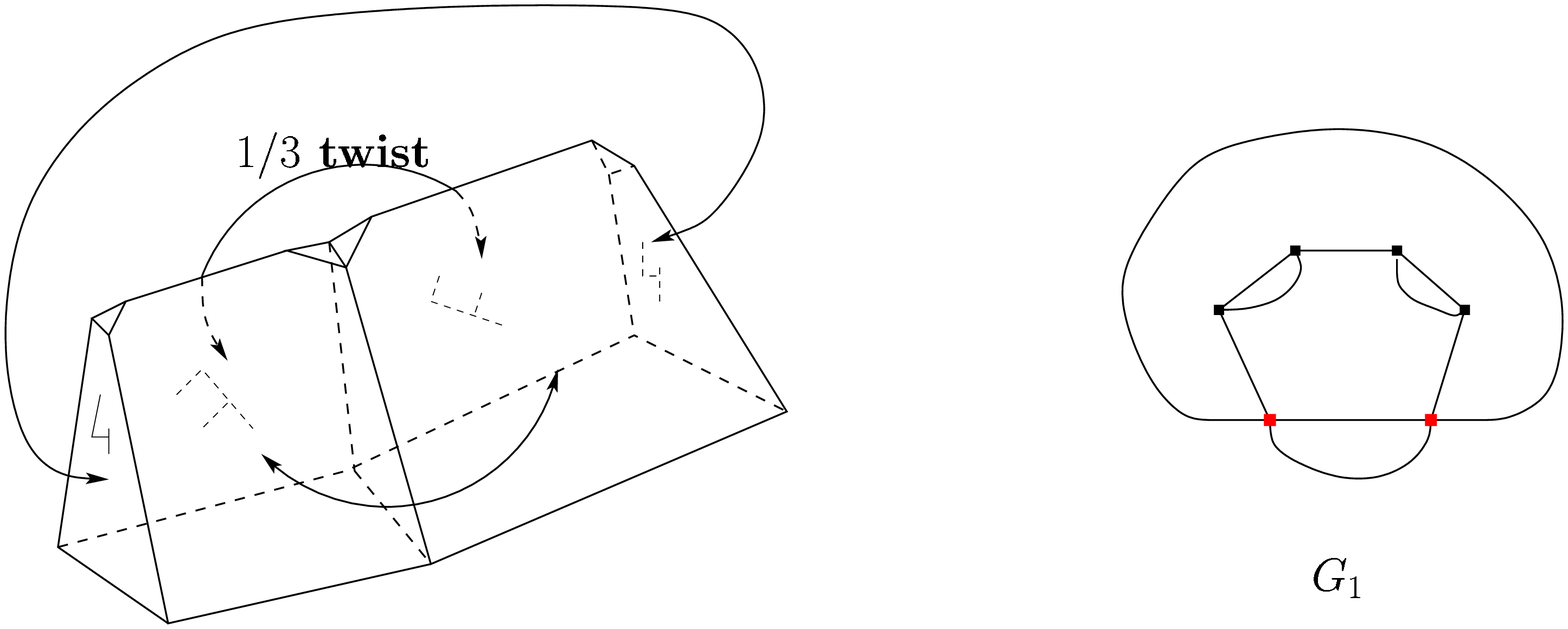, width= 2.9in} \caption{A
cell decomposition of a punctured $L(3,1)$ and the dual $G_1$ of
its boundary.} \label{fig4.29}}
\end{figure}

Since $G'$ is planar, it is the dual of a cell decomposition of a
2-sphere, and hence, we did not change the topology of $M'$.
Consider the cell decomposition $N'$ of the punctured $L(3,1)$ in
figure~\ref{fig4.29} with the graph $G_1$ corresponding to the
dual of $\partial N'$. We cut along the face labeled 4 in figure~\ref{fig4.28},
and we insert truncated tips so that the resulting graph
 $G'_1$ of $G_1$ is
colored homeomorphic to $G'$. It is now clear that the union of
$M'$ and $N'$ gives a well-defined triangulation of $M \# L(3,1)$
with $(n+2)$ tetrahedra.

\begin{construction} \label{cons7}
Let M be a triangulated closed orientable 3-manifold with $|M| = n$ and at least two vertices.
Then there exists a triangulation for $M \# \mathbf{RP}^{3}$ with (n + 2) tetrahedra.
\end{construction}

\bigskip
\begin{figure}[h]
\centering{\psfig{figure=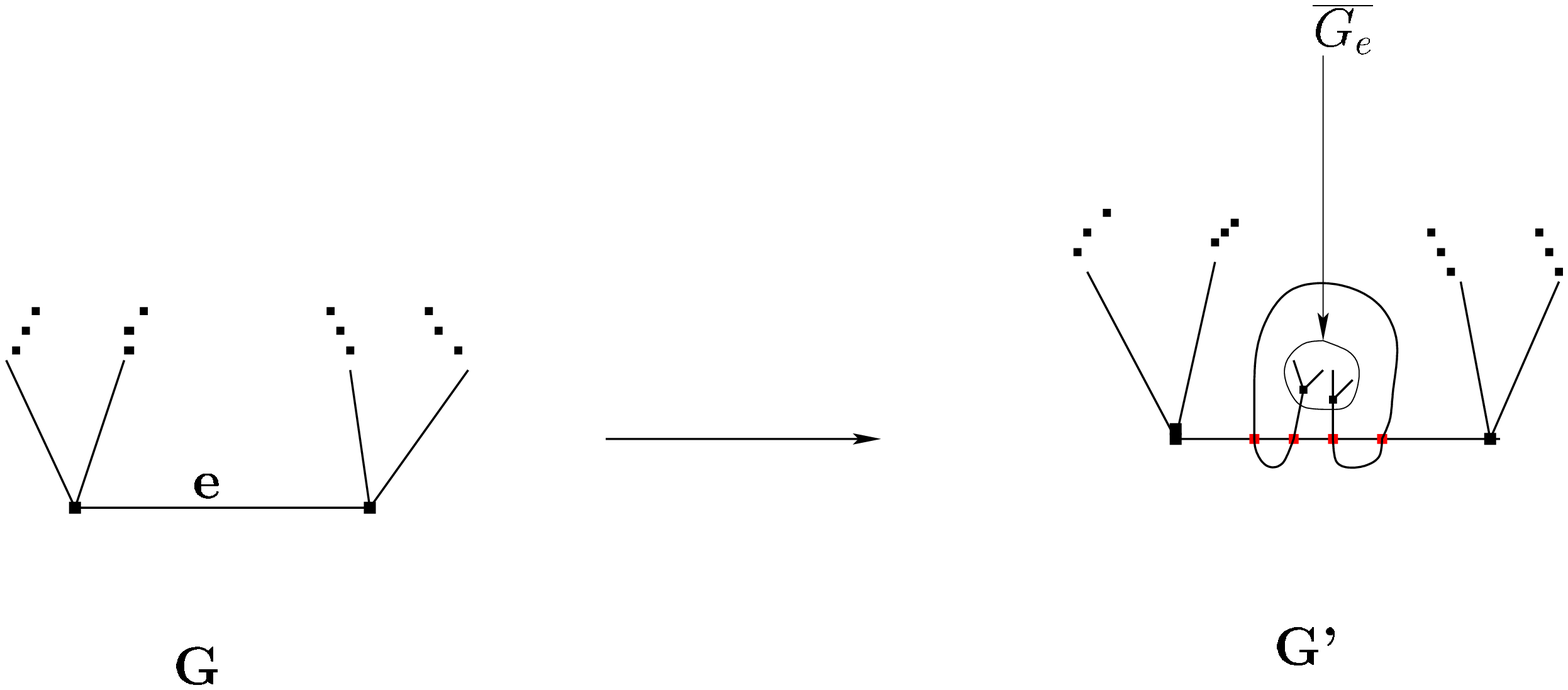, width= 3in} \caption{The new triangulation $G'$ of $G$.}\label{fig4.25}}
\end{figure}

\bigskip
\begin{figure} [h]
\centering{\psfig{figure=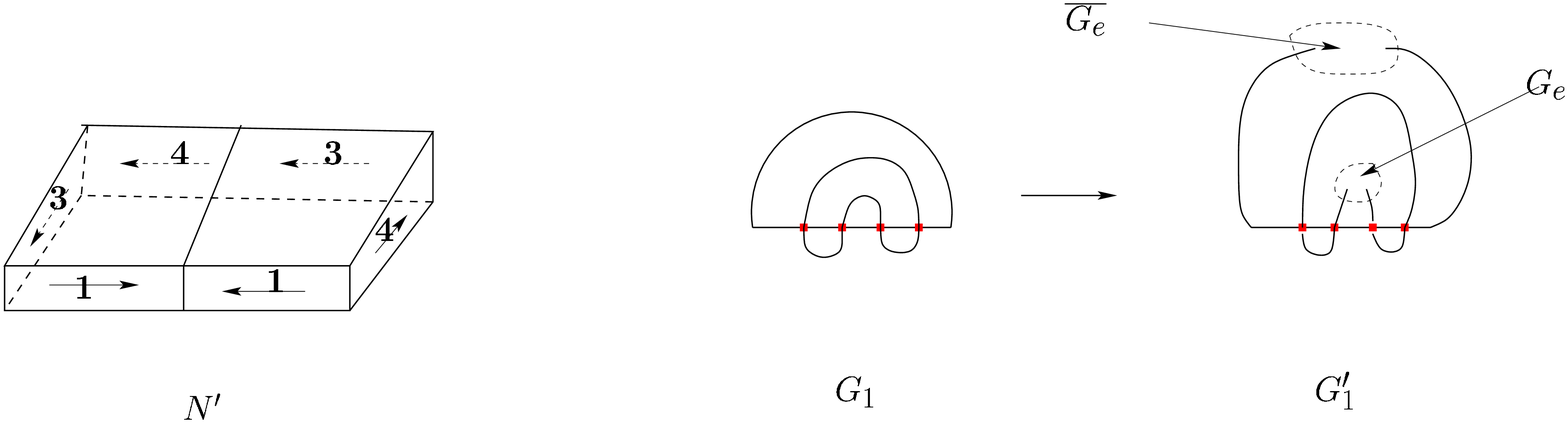, width= 4.3in} \caption{A
cell decomposition of a punctured \R and the dual $G_1$ of
 its boundary.}\label{fig4.26}}
\end{figure}

 Let $M'$ be the manifold obtained from $M$ after removing the link of a vertex. Without loss of generality, we assume
 we can find a truncated tetrahedron as in figure~\ref{fig4.1}. Let $T$ be the triangulation
 of $\partial M'$, and $G$ its dual. Let $e$ be the thickened edge from Figure~\ref{fig4.1}. We insert 4 prisms
 along the shaded face and we insert tips of tetrahedra so that the
 resulting graph $G'$ of $G$ correspond to figure~\ref{fig4.25} above. Note, $G'$ is planar and so we did not
  changed the topology of $M'$. Consider now the cell
 decomposition $N'$ of the punctured $\mathbf{RP}^3$ in figure~\ref{fig4.26}.

We cut along the face labeled 4 in figure~\ref{fig4.26}, and we insert truncated tips so that the
resulting graph $G'_1$ of $G_1$ is colored homeomorphic to $G'$.
It is now clear that the union of $M'$ and $N'$ gives a
well-defined triangulation of $M \# \mathbf{RP}^3$ with $(n+2)$
tetrahedra.

\begin{construction} ~\label{cons8}
Let M be a triangulated closed orientable 3-manifold with $|M| = n$ and at least three vertices.
Then there exists a triangulation for $M \# (S^{1} \times S^{2}$) with (n + 2) tetrahedra.
If M contains a good vertex, this construction can be done with one extra tetrahedron only.
\end{construction}

\begin{prop}
There exists a 1-vertex triangulation of any closed orientable 3-manifold.
\end{prop}

\underline{{\bf Proof :}} Let $M$ be a closed orientable reducible 3-manifold equipped with a
 triangulation. Suppose $M$ has more than 1 vertex. Consider the 1-vertex 1-tetrahedron triangulation of
$S \sp 3$. Then $M \# S^{3}$ has $(t + 3)$ tetrahedra and $(v-1)$
vertices by our Construction~\ref{cons1}. We can repeat this
construction $(v-2)$ times and obtain a manifold homeomorphic to
$M$ with $(t+3(v-1))$ tetrahedra and 1-vertex only. Q.E.D.

\section{Small normal 2-Spheres in Minimal Triangulations}

 Using Theorem~\ref{MT} and the constructions in the previous section, we show that minimal triangulations of
 reducible 3-manifolds contain {\it small} non-trivial normal 2-spheres. By \textbf{\textit{small}}, we mean
  normal 2-spheres which have quadrilateral types in not more than 2 tetrahedra. This will be essential in Andrew
  Casson's algorithm to check if a minimal triangulation is reducible or not.

\smallskip

\begin{defn} A 3-manifold $M$ is said to have a \textbf{\textit{minimal triangulation}} $\tau$ if $\tau$
 contains the smallest number of tetrahedra over all possible triangulations of $M$. By abuse of language, we
 say that $M$ is minimal.

\smallskip


\smallskip

 We say that a normal surface $S$ has $n$ \textbf{\textit{quadrilateral types}} ( $<S>$ = n) if
there exist exactly n tetrahedra with the property that S intersects each of these tetrahedra in quadrilaterals.
Note, if S and T are two compatible normal surfaces with $<$S$>$ = n and $<$T$>$ = m, then
 $(n+m) \geq$ $<$$S + T$$>$ $\geq $max(n, m). Moreover, normalizing  an embedded surface decreases its weight but
  may increase its number of quadrilateral types. Given a normal surface F with normal vector $x_F$, let $q_1$,
  ..., $q_{4t}$ be the entries in $x_F$ corresponding to the quadrilaterals of F.
  We denote by $\mathbf{\#_q (F)}$, the sum of the $q_i$'s.

\end{defn}

To prove the existence of small normal 2-spheres in minimal triangulations, we first
 show the existence of a normal 2-sphere with the property that all the collapsing surfaces are inessential.
  If we find such a 2-sphere, then collapsing it will result in a
  decomposition of M in exactly two summands. Because $M$ is minimal, we will use the constructions in the
  previous section to conclude that $S$ cannot have more than 2 quadrilateral types. To show the existence of a
  non-trivial normal 2-sphere with no essential collapsing surfaces, we
 use a result due to W. Jaco and J. Tollefson ~\cite{JT}.

\smallskip

 \underline{Theorem 4.1} (~\cite{JT}): A normal two-sphere $F$ is a vertex surface if and
only if $F$ has the property that whenever there exists an annulus $A$ which is an exchange surface for $F$
 then the two disjoint disks in $F$ bounded by $\partial A$ are normal isotopic.

 Here, an \textbf{\textit{exchange surface}} $A$ for $F$ is a surface with the following
 properties: 1) $fr(A) = A \cap F$, 2) $A$ has an orientable regular neighborhood $N(A)$, and 3) for
  every tetrahedron $\Delta$, each component of $\Delta \cap A$ is a 0-weight disk $L$ spanning two
  distinct elementary disks $E_1$, $E_2$ of $F$ such that $\partial L = L \cap (E_1 \cup E_2 \cup \partial \Delta)$
  and $L \cap E_i$ is an arc joining the interiors of two distinct 2-faces of $\Delta$.

 We see from this definition that if $A$ is a collapsing annulus or a collapsing Mobius band,
then we can either push $A$ off the 2-skeleton (if $A$ is an
annulus) or look at $\partial Nbhd(A)$ (if $A$ is a Mobius band)
to find an annulus which is an exchange surface. Hence, if $F$ is
a vertex 2-sphere, we now know that any collapsing annulus or Mobius
band
 is inessential. On the other hand, Jaco and Tollefson's theorem does not say anything about collapsing disks.
 To go around this problem, we define a complexity $Q$. Let $S$ be a normal surface. Then \textbf{\textit{Q(S)}}
 represents a pair, ordered lexicographically, whose first and second entries are $< S >$ and $\#_q (S)$
 respectively: Q(S) = ($<S>, \#_q (S)$).

We are now ready to prove the Lemma for the main theorem.

\begin{lem} \label{L4}
Let M be a triangulated closed orientable 3-manifold which contains a non-trivial
normal 2-sphere. Then M contains a non-trivial normal 2-sphere whose collapsing surfaces are all inessential.
\end{lem}

 \begin{proof}
 Consider the set of non-trivial normal 2-spheres.
This set is non-empty by assumption. In this set, choose the 2-sphere $F$ which is minimal
 with respect to $Q$. We want to show that such a 2-sphere does not have any essential collapsing
 surfaces, but first, we want to show that it is a vertex surface.

 Suppose $F$ is not a vertex surface, i.e. suppose $x_F$ is not a vertex of \textit{P(M, T)},
i.e. suppose for all positive integer $k$, $k \cdot x_F$ is not a vertex solution. Let $k$ be the smallest
 positive integer such that $k \cdot x_F = S$ is an integral solution (and, by assumption, not a vertex
 solution). Then $\chi(S)$ must divide $\chi(F)$, and so $\chi(S) = 1$ or 2.

 \underline{Case 1:} $\chi(S) = 2$. Since $F$ is an integral multiple of $S$, $S = F$. Suppose
there exists a positive integer $n$ such that $nF = V_1 +$... $+ V_k + W_1 +$... $+W_r + X$, where at least
 one of the summands is not an integral multiple of $F$. Without loss of generality, we can assume that the
 $V'_i$s are 2-spheres, the $W'_i$s are real projective planes, and $X$ is a (possibly empty or disconnected) surface
 with non-positive Euler characteristic. Because the Euler characteristic is preserved
  under surface addition, we have $2k + r \geq 2n$.

 \underline{Subcase 1:} Suppose that $r = 0$. Then $nF = V_1 +$... $+ V_k + X$ with $k \geq n$.
 Consider the 2-sphere, say $V_1$, which has the smallest $\#_q$ over all the $V_i$'. Since the surface addition
  preserves the quadrilateral types and the number of quadrilaterals for each quadrilateral type, we have $< V_1 >
  \leq < F >$. But by assumption, $< F > \leq < V_1>$. Hence, $< F > = < V_1 >$. It is crucial to notice that, not
  only $F$ and $V_1$ have the same number of quadrilateral types, but they must also have their quadrilaterals in the same
   tetrahedra. Moreover, $n \cdot (\#_q ( F )) \geq k \cdot (\#_q ( V_1 )) + \#_q ( X )$. Since $k \geq n$ and
    $\#_q (V_1) \geq \#_q(F)$, the only way for the inequality to be true is if it is an equality and $k = n$
    and $\#_q( X ) = 0$. We conclude that $X$ was actually empty. Also, $\#_q (V_1) \leq \#_q(F)$. By assumption,
     $\#_q (F) \leq \#_q(V_1)$, so $\#_q(F) = \#_q(V_1)$. Therefore, $V_1$ is a parallel copy of $F$. We remove
      $V_1$, X,  and a copy of $F$ from the equation to obtain $(n-1)F = V_2 +$... $+ V_n$. We repeat the argument
      for the next 2-sphere, say $V_2$, which as the smallest $\#_q$. We conclude that $V_i$ is a parallel copy of
       $F$ for $1 \leq i \leq n$. This contradicts the fact that $F$ was not an integral solution.

\underline{Subcase 2:} Suppose $ r \neq 0$. We look at the equation
$2nF = 2V_1 +$... $+ 2V_k + 2W_1 +$... $+2W_r + 2X$, where $2V_i$ represents two copies of a 2-sphere and $2W_j$
 represents a 2-sphere. We repeat the same argument as in Subcase 1 to conclude that each $V_i$ is a parallel copy
 of $F$, that each $2W_i$ is a parallel copy of $F$, and that $X$ is empty. Since, say $2W_1$ is parallel to $F$,
 we conclude that $F$ is the double of the projective plane $W_1$. This contradicts the fact that $k$ was the smallest
  positive integer such that $S$ (= $F$) is an integral solution.

\smallskip

\underline{Case 2:} $\chi(S) = 1$. Since $F$ is an integral multiple of $S$, $F = 2S$. Suppose
there exists a positive integer such that $nS = V_1 +$... $+ V_k + W_1 +$... $+W_r + X$, where the $V_i$'s are
 2-spheres, the $W_j$'s are real projective planes, and $\chi(X) \leq 0$. We look at the new equation
 $nF = 2nS = 2V_1 +$... $+ 2V_k + 2W_1 +$... $+ 2W_r + 2X$. We run the same argument as in Subcase 1 of Case 1 to
  conclude that each $V_i$ and each $2W_j$ is a parallel copy of $F$, and that $X$ is empty. Hence, each $V_i$ and
  $2W_j$ represents the double of $S$. Therefore, each summand of the equation is an integral multiple of $S$. We
  conclude that $S$ is a vertex solution, or equivalently, $F$ is a vertex surface. Contradiction.

\medskip

 Using the above theorem from Jaco and Tollefson, we know that if $F$ is a non-trivial normal
2-sphere which minimizes $Q$, then $F$ does not have any essential
collapsing annuli or Mobius bands.

\smallskip

 Let $D$ be an essential collapsing disk. We construct another non-trivial normal 2-sphere with a smaller number
of quadrilaterals: let $D_1$ and $D_2$ be the two disjoint disks
on $F$ such that $D_1 \cap D_2 = \partial D$. Consider the
2-sphere $D \cup D_1$. By assumption, this 2-sphere does not bound
a ball whose interior is disjoint from $F$. We push $D$ off the
vertex and we obtain a new normal 2-sphere $F'$. See figure~\ref{fig4.20} below. Note, $F'$ is
represented by the same set of quadrilaterals as $D_1$. Because
$\#_q(F)$ is minimal, we conclude that  $D_2$ is composed of
triangles only. Hence, $D \cup D_2$ bounds a ball. This is a
contradiction and so $D$ was inessential after all.

\bigskip
\begin{figure}[h]
\centering{\psfig{figure=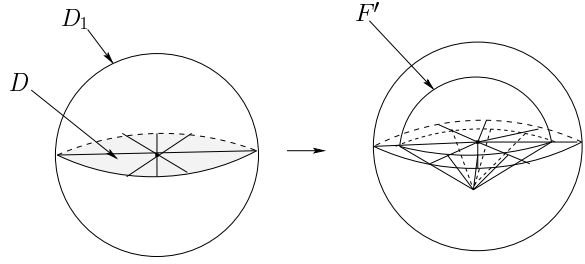, width= 3in} \caption{$D$ is pushed off the 2-skeleton.}\label{fig4.20}}
\end{figure}
\vspace{.2in}

\end{proof}

\begin{rem}
 Let $A$ be a collapsing surface. We now know that the two disjoint disks $D_1$ and $D_2$
on $F$ such that $\partial A = \partial D_1 \cup \partial D_2$ are normal isotopic. Moreover, Lemma 4.6 of
~\cite{JT} tells us that the ball with boundary $A \cup D_1 \cup D_2$ whose interior is disjoint from $F$ does not
contain any vertex from the triangulation. This is a crucial part of the proof of the following theorem.

\end{rem}

\bigskip


\begin{thm} ~\label{MT3}
Let M be a closed orientable 3-manifold equipped with a minimal
triangulation. If M contains a non-trivial normal 2-sphere, then M
contains such a 2-sphere with at most 2 quadrilateral types.
\end{thm}

By Lemma~\ref{L4} we know there exists a non-trivial normal
2-sphere $S$ with only inessential
  collapsing surfaces. If $<S> \leq 2$, then we are done. So we will assume that $<S>$ $\geq 3$ and we will
   contradict
   the minimality of $M$. Consider the 2-sphere $F$ constructed in Lemma~\ref{L4}.
   We collapse $S$. Since $S$ has no essential collapsing surfaces, we end up with a decomposition of $M$ into
   exactly 2 summands:

\begin{enumerate}

 \item $M \cong M_{1} \# M_{2}$.

 \item $M \cong M_{1} \# (S^{1} \times S^{2})$

\item $ M \cong M_1 \# L(3, 1)$.

\item $M \cong M_1 \# \mathbf{R}P^3$.

\item $M \cong M_1 \# S^3$.

\item $M \cong L(3, 1) \# L(3, 1).$

\item $M \cong L(3, 1) \# \mathbf{R}P^3.$

\item $ M \cong L(3, 1) \# S^3.$

\item $M \cong \mathbf{R}P^3 \# \mathbf{R}P^3.$

\item $M \cong \mathbf{R}P^3 \# S^3.$

\item $M \cong S^3 \# S^3.$

\item $M \cong L(3, 1) \# (S^{1} \times S^{2})$.

\item $M \cong \mathbf{R}P^3 \# (S^{1}\times S^{2}).$

\item $M \cong S^3 \# (S^{1} \times S^{2}).$
 \end{enumerate}

\smallskip

By the above remark, the sum of the number of vertices in the resulting triangulated summands of $M$ must be at
least the number of vertices in the original triangulation. Hence, we can eliminate (6), (7), (9), (12), and (13).

Jaco and Rubinstein ~\cite{JR} showed that the minimal triangulations of $S^3$, $\mathbf{R}P^3$, $S^1 \times S^2$,
and $L(3, 1)$ have less than 3 tetrahedra. Hence, any non-trivial normal
2-sphere in these triangulations has less than 3 quadrilateral types.
 We can thus eliminate (8),
(10), (11), (14).

(5) clearly contradicts the minimality of $M$.

In (3) and (4), $M_1$ must contain at least the same number of
vertices as $M$ plus one from the collapsing of $S$.  Hence, $M_1$
must have at least 2 vertices, and constructions ~\ref{cons6} and
~\ref{cons7} can be used to obtain a new triangulation of $M$ with
less tetrahedra. Contradiction.

For similar reasons, constructions ~\ref{cons1} and ~\ref{cons8}
can be used in (1) and (2) to construct triangulations of $M$ with less
tetrahedra. Contradiction. This proves the theorem. \qed

\section{Casson's Algorithm to Decompose a Closed Orientable 3-Manifold into Irreducible Pieces}

In 1952, Edwin E. Moise ~\cite{Mo} proved that any 3-manifold can be triangulated and that,
 given any two triangulations $K_{1}$ and $K_{2}$ of the same 3-manifold, $K_{1}$ and $K_{2}$ are equivalent, i.e.
  they have isomorphic subdivisions. As of today, it has actually been proven that we can obtain $K_{2}$ from $K_{1}$
  through a series of Pachner moves ~\cite{Pa}, though no bounds on the number of moves has been found for
   arbitrary 3-manifolds \footnote{Explicit bounds have been found for Seifert fibered spaces in ~\cite{Mi}}.

Let $M$ be a closed orientable 3-manifold triangulated with $t$ tetrahedra. We mention that, even
 though the following lemmas were proved by the author, the actual results are due to Andrew Casson.

\bigskip

\begin{defn}

Let \emph{$\overline{N}$}(M) be the set of surfaces in M which intersect each tetrahedron in triangles or
 quadrilaterals. Precisely, \emph{$\overline{N}$}(M) is the set of surfaces which satisfy the normal surface
  equations but may not satisfy the quadrilateral property.

A type w is a function which assigns, to each tetrahedron, one of
the 3 possible types of quadrilateral. Let N(M) be the set of
normal surfaces and let $N_{w}$(M) be the set of normal surfaces
of type w. Note  that $N(M)=  \bigcup_{w}  N_{w}(M)$.
\end{defn}

 The reason we look at $\bigcup_{w} N_{w}(M)$ instead of $N(M)$ is because the solution
 space of $N(M)$ does \underline{not} form a cone in $\mathbf{R}^{7t}$ whereas each $N_{w}(M)$ does
 form a cone. Indeed, if $S_{1}$, $S_{2} \in N(M)$ are not of the same type, their normal sum
 do not represent an embedded surface.

\bigskip
\begin{claim} ~\label{C6.3}
There exists a 2-sphere in $N(M)$ if and only if there exists a
surface $S$ in  $N_{w}(M)$ for some $w$, with $\chi(S)>0$.
\end{claim}

\underline{\bf{Proof:}} $\Rightarrow$ Trivial since $N(M)=
\bigcup_{w}  N_{w}(M)$.

\hspace{.35in} \space $\Leftarrow$ \space If there is a surface
$S$ in $N_{w}(M)$ for some $w$, then this surface must be
 embedded. Since $\chi(S)>0$, one of the connected components of $S$ must be homeomorphic to either a $\mathbf{RP}^{2}$
  or a 2-sphere. If it is a 2-sphere, then  we are done. If it is a $\mathbf{RP}^{2}$, then a
  regular neighborhood of this projective plane is homeomorphic to a twisted $I$-bundle over it. This twisted $I$-bundle
   is homeomorphic to a punctured $\mathbf{RP}^{3}$. Hence, $\partial N(RP^{2})\cong S^{2}$ which has the desired
   property. \qed

Suppose M has a unique vertex $v$. Let $t_{1}, t_{2}, ..., t_{4t}$ be the triangle coordinates corresponding to $v$.
Let $N_{w, i}(M)=\{ S \in N_{w}(M) | t_{i}=0 \}$

\bigskip

\begin{claim} ~\label{C6.4}
There exists a non-trivial 2-sphere in $N_w(M)$ if and only if
there exists a surface $S$ in $N_{w, i}(M)$ for some w, $i$, and
with $\chi (S)>0$.
\end{claim}

\underline{\bf{Proof:}} $\Rightarrow$ If $S^{2}$ is non-trivial, then it cannot be the boundary of a neighborhood
 of a vertex and hence, we must have $t_{i}=0$ for some $i$.

\hspace{.35in} \space $\Leftarrow$ if $S \in N_{w, i}(M)$ then $S$ cannot be the boundary of a
 vertex neighborhood and hence, cannot be trivial. As in claim~\ref{C6.3}, $\chi (S)>0$ implies the existence of a
 non-trivial $S^{2}$. \qed

\bigskip
 We define $\mathbf{C_{w, i}(M)}$ to be the solution space to the surface equations with type $w$, $t_{i}=0$ for a
  fixed $i$,
 and the rest of the variables being \emph{real} non-negative.

\underline{Remark:} let us define the Euler characteristic
function on a rational-entry vector satisfying the surface equations and the quadrilateral property. Let
$S=(a_{1}, a_{2}$, ..., $a_{7t})$ be an integer solution of the surface equations representing an embedded surface.
Let $a_{1}$, ..., $a_{4t}$ and $a_{4t+1}$, ..., $a_{7t}$ denote the respective coefficients of the triangle and
quadrilateral entries. We then define $\chi$ in the following way:

$\displaystyle{\chi (S) = \lbrack \sum_{i=1}^{4t}a_{i} \cdot (\frac{\rm 1}{\rm d_{i_{1}}} + \frac{\rm 1}{\rm d_{i_{2}}} +
 \frac{\rm 1}{\rm d_{i_{3}}}) - \sum_{i=1}^{4t}a_{i} \cdot (\frac{\rm 3}{\rm 2}) + \sum_{i=1}^{4t}a_{i} \rbrack}$ $\displaystyle{+}$

 \hspace{1cm} $\displaystyle{\lbrack \sum_{i=4t+1}^{7t}a_{i} \cdot (\frac{\rm 1}{\rm e_{i_{1}}} + \frac{\rm 1}{\rm e_{i_{2}}} + \frac{\rm 1}{\rm
   e_{i_{3}}} + \frac{\rm 1}{\rm e_{i_{4}}}) - \sum_{i=4t+1}^{7t}a_{i} \cdot (\frac{\rm 4}{\rm 2}) + \sum_{i=4t+1}^{7t}a_{i}
    \rbrack}$

     $\displaystyle{\chi(S) = \lbrack \sum_{i=1}^{4t}a_{i} \cdot (\frac{\rm 1}{\rm d_{i_{1}}} + \frac{\rm 1}{\rm d_{i_{2}}} + \frac{\rm 1}
    {\rm d_{i_{3}}} - \frac{\rm 1}{\rm 2}) \rbrack + \lbrack \sum_{i=4t+1}^{7t}a_{i} \cdot (\frac{\rm 1}{\rm e_{i_{1}}} +
    \frac{\rm 1}{\rm e_{i_{2}}} + \frac{\rm 1}{\rm e_{i_{3}}} + \frac{\rm 1}{\rm e_{i_{4}}} - 1) \rbrack}$, where the
    $d_{i}$'s and the $e_{j}$'s represent the respective degrees of the edges touching the triangle $i$ and the quadrilateral $j$.

\hspace{1cm} Now let $V$ be a rational solution of the surface equations, satisfying also the quadrilateral property
 ( $V$= ($\frac{\rm a_{1}}{\rm b_{1}}, \frac{\rm a_{2}}{\rm b_{2}}$, ..., $\frac{\rm a_{7t}}{\rm b_{7t}}) $).
 Let $N$= $\Pi_{i=1}^{7t}b_{i}$. We then define $\chi$ as follows: $\chi (V) = \frac{\rm 1}{\rm N}\chi ( N \cdot V)$,
  where $N \cdot V$ denotes the integer solution $(\frac{\rm a_{1} \cdot N}{\rm b_{1}}, \frac{\rm a_{2} \cdot N}{\rm b_{2}}$,
   ..., $\frac{\rm a_{7t} \cdot N}{\rm b_{7t}})$. From this definition, it follows that $\chi$ is a linear function.
   Indeed, let $V_{1}$ and $V_{2}$ be two rational solutions. Then, $\chi (V_{1}) + \chi (V_{2}) = \frac{\rm 1}{\rm MN}
    \chi (\frac{\rm 1}{\rm MN}(V_{1} + V_{2})) = \frac{\rm 1}{\rm MN} \chi (MN\cdot V_{1}) + \frac{\rm 1}{\rm MN}
     \chi (MN\cdot V_{2}) = \frac{\rm 1}{\rm M}  \chi (M\cdot V_{1}) + \frac{\rm 1}{\rm N} \chi (N\cdot V_{2}) =
     \chi (V_{1}) + \chi (V_{2})$.

\bigskip
\begin{claim} ~\label{C6.5}
There exists a surface S in $N_{w, i}(M)$ for some $w, i$ with
$\chi (S)>0$ if and only if there exists a vector $V$ in $C_{w,
i}(M)$
 for some $w , i $, such that $\chi (V)>0$, $V$ is on an external ray of the cone $C_{w,i}(M)$, and $V \in \mathbf{Z}^{7t}$.
\end{claim}

\underline{\bf{Proof :}} $\Leftarrow$ Suppose $\exists V \in C_{w, i}(M)$, $\chi (V)>0$, and $V \in \mathbf{Z}^{7t}$.
Then , trivially, $V \in N_{w, i}(M)$.

\hspace{.35in} \space $\Rightarrow$ If $S \in N_{w, i}(M)$, then
$S \in C_{w, i}(M)$. Suppose now that $S$ is not on
 some external ray, i.e. suppose that $S= V_{1} + V_{2} + ... + V_{n}$ where $V_{j} \in C_{w, j}(M)$, and the $V_{j}$'s
 are connected and they lie on an external ray. Since $S$ has integer entries, we can assume without
 loss of geneality
 that the $V_{j}$'s have rational entries.
 Then $\chi (S) = \chi (V_{1} + ... + V_{n})= \chi (V_{1}) +... + \chi (V_{n}) >0$. This implies that $\chi (V_{j}) >0$
  for some $j$. We multiply each entry of $V_{j}$ be the least common multiple of the denominator of the entries.
   Since $\chi (V_{j}) >0$, we conclude that either $V_{j}$ is homeomorphic to a 2-sphere or a projective plane.
    If it is a projective plane, then we double each entry in the vector representation of the surface to
    obtain a 2-sphere on the same external ray. Also, $V_{j}$ cannot represent a trivial 2-sphere because $i = 0$
     for some $i$. $V_{j}$ or its double gives us the desired surface. $\qed$

\begin{defn}
$S \subset M$ is called \emph{almost $^{2}$ normal} (resp. $almost$ $normal$) if it is normal and if there exists
at most one tetrahedron $\tau$ in which S intersects $\tau$ in triangles and octagons (resp. one octagon) only.

 Define $A_{w, i, l}(M)$ to be the set of surfaces of type $w$, with $t_{i}=0$, and one type of octagon
being allowed only in the $l^{th}$ tetrahedron.
\end{defn}

\medskip
\begin{claim} ~\label{C6.7}
$\exists$ an almost normal $S^{2}$ $\Leftrightarrow$ $\exists S \in A_{w, i, l}(M)$ for some $w$, $i, l$ with
$\chi (S) - o(S) >0$, where o(S) is the number of octagons in S in the $l^{th}$ tetrahedron.
\end{claim}

 \underline{\bf{Proof :}} $\Rightarrow$ If there is an almost normal $S^{2}$, then clearly
$S \in A_{w, i, l}(M)$ for some $w, i, l$. Since $S$ is almost normal, then $o(S)=1$. This implies that
$\chi (S) - o(S)=1 >0$.

\hspace{.4in} \space $\Leftarrow$ If there is a $S$ in $A_{w, i, l}(M)$ with $\chi (S)- o(S) >0$, then, by definition,
$S$ is almost$^{2}$ normal. Since $\chi (S) - o(S) >0$, and o(S)$\geq 1$ (if o(S)=0 then $S$ is normal non-trivial, but
 we already have cut along all such surfaces), then $\chi (S) > 1$. If $\chi (S) = 2$, then $o(S) = 1$ and we have an
 almost normal $S^{2}$. If $\chi (S) >2$, then $S$ must have an $S^{2}$ component with $\chi (S^{2}) - o(S^{2}) > 0$.
 This again implies that $o(S^{2}) = 1$ and we have an almost normal $S^{2}$.\qed

\bigskip

\underline{\underline {\bf{The algorithm :}}}

Let $M$ be a closed orientable 3-manifold equipped with a $t$-tetrahedra one-vertex triangulation. The set of
surface equations is completely determined by the triangulation, and so are the cones $C_{w}(M)$. There
 are  $3^{t}$ of them, one for each type $w$.

 The first task of the procedure is to find a non-trivial normal 2-sphere. We use claim~\ref{C6.5} to find one.
 Given a cone $C_{w, i}(M)$, we look at the set $A= C_{w, i}(M)$ $\bigcap$
    $\{ \sum_{i=1}^{7t} t_{i} = 1  \} $. This set is a convex compact polyhedron with
     vertices having rational entries. On this polyhedron, we maximize the Euler characteristic
   function defined above.  Linear programming theory tells us that this function attains its maximum at a vertex of this
    polyhedron and hence, on an external ray of the cone $C_{w, i}(M)$. There are several methods to find such a maximum.
    Indeed, Schrijver ~\cite{Sch} (Theorem 15.3 page 198) first proved the existence of a method to find the
    maximum of a linear function on a convex compact polyhedron with a running time polynomial in the size of a matrix $B$.
    In our context, $B$ describes the surface equations and the non-negativity of the entries $t_{i}$. The size of $B$
     (as defined in Schrijver page 29) happens to be a polynomial with respect to $t$. In fact, $size(B) \leq 182t^{2}$.
     This upper bound comes from the fact that B is a $(7t) \times (6t + 7t)$ matrix with integer entries smaller than 2.
      See ~\cite{Kar}
      for an explicit algorithm.

      This vertex, on which $\chi$ is
      maximum, has rational entries and so we multiply it by the least common multiple of the denominator of the entries.
      This new vector $S$ with integer entries can have Euler characteristic 0, 1, or 2:

$\chi(S)=1$: $S$ represents a projective plane. We look at the surface $2S$ which must represent a non-trivial 2-sphere
 since $S$ is non-trivial.

$\chi (S)=2$: $S$ is a non-trivial 2-sphere.

If $\chi (S)<0$, then by claim~\ref{C6.4} and claim~\ref{C6.5}
there are no non-trivial $S^{2}$ and the procedure stops here.

$\chi (S)>2$  would contradict the facts that $S$ is on an
external ray of $C_{w, i}(M)$.

Once the procedure finds a non-trivial 2-sphere ($S$ or $2S$), it
collapses it using theorem~\ref{MT}.  After cutting and collapsing
along all non-trivial $S^{2}$'s, the procedure needs to perform
one more task. It needs
 to check if some of the resulting summands are homeomorphic to $S^{3}$.

 The second task of the procedure is to find an almost normal 2-sphere. To do that, we use the Thompson-Rubinstein theorem
  which states that one of the resulting pieces is homeomorphic to $S^{3}$ if and only if there exists an almost normal
   2-sphere in it. To find such a 2-sphere, we maximize the linear function $\chi (\bullet) - o(\bullet)$
   over the set $A_{w, i, l}(M) \bigcap$
   $\{ \sum_{i=1}^{7t} t_{i} = 1 \}$. This gives us a solution on a rational vertex $W$, which in turn gives us an
   integer vector $V$. If we obtain $\chi (V) - o(V) \leq 0$, then by claim~\ref{C6.7} there are no almost normal
   2-spheres and $M_{i}$ is not homeomorphic to $S^{3}$ and the procedure stops.  If $\chi (V) - o(V) = 1$,
    then $M_{i}$ is homeomorphic to $S^{3}$. Note, $\chi (V) - o(V) > 1$ would imply that
   $\chi(V) > 2$ which would mean that $V$ is disconnected which is impossible since $V$ is a vertex.

Before calculating the complexity of the algorithm, we need to show that the procedure terminates after finitely many
steps.

\underline{\bf{Fact 1:}} The procedure terminates after cutting along and collapsing finitely many non-trivial
2-spheres. This follows from Theorem~\ref{MT} since cutting along a non-trivial $S^{2}$ strictly reduces the
 original number of tetrahedra in $M$.

\underline{\bf{Fact 2:}} If there are no non-trivial $S^{2}$, the
procedure terminates.  Indeed, claim~\ref{C6.5} tell us that there
are no such 2-sphere if and only if $\chi (S) \leq 0$.

\underline{\bf{Fact 3:}} The procedure cannot be looking indefinitely for an almost normal $S^{2}$ if there are
 none. Indeed, if $\chi (V) - o(V)\leq 0$, then by claim~\ref{C6.7} there are no almost normal 2-spheres.

 We summarize the steps needed in the algorithm to decompose a closed orientable 3-manifold into
irreducible pieces. Let $M$ be given by a $t$-tetrahedra and $v$-vertex triangulation.

 \underline{\bf{Step 1:}} Construct a (possibly disconnected) normal surface $S$ obtained by
normalizing the boundary of a regular neighborhood of a maximal tree of the 1-skeleton. We collapse $S$ using
Theorem~\ref{MT}. This may result in a decomposition of $M$: $M$  $\cong  M_{1}$ $\# M_{2}$...$\# M_{k}$
 $\# r_{1}(S^{1}\times S^{2})$ $\# r_{2}\mathbf {RP}^{3}$ $\#r_{3}L(3,1)$. For each summand not homeomorphic to
 $\mathbf{RP}^{3}$, $S^{1} \times S^{2}$, $S^{3}$, or $L(3,1)$, and having more than one vertex in its
  triangulation repeat step 1. We end up with a decomposition of $M$ where each
summand is either $\mathbf{RP}^{3}$, $S^{1} \times S^{2}$,
  $L(3,1)$, or it has a 1-vertex triangulation. Let us call $M$ one of the 1-vertex summands. We will repeat the procedure
   for each of the other 1-vertex summands.

\underline{\bf{Step 2 :}} Go through each cone $C_{w,i}(M)$ to find a non-trivial 2-sphere. If no such sphere
 is found, go to step 3. If one is found, then collapse it using Theorem~\ref{MT}. We obtain a further
  decomposition of $M$. For each new summand having more than 1-vertex in its triangulation, repeat step 1.
For each new summand having a 1-vertex triangulation, repeat step 2.

\underline{\bf{Step 3:}} Go through each cone $A_{w,i,l}(M)$ to find an almost$^2$ normal 2-sphere.
  If none are found, go to step 4. If one is found, then we know the summand is homeomorphic to $S^{3}$.

\underline{\bf{Step 4:}} We obtain a decomposition of $M$ where each summand is not homeomorphic to $S^3$  and is either
irreducible with a 1-vertex triangulation, or is homeomorphic to $\mathbf{RP}^{3}$, $S^{1} \times S^{2}$, or $L(3,1)$.

We now describe the complexity of the algorithm. Note, our goal here is not to give explicit bounds for the
running time of this algorithm. Rather, we only make a distinction between polynomial and exponential
running time.

\underline{\underline{\bf{Complexity of the algorithm:}}}

\noindent
\begin{enumerate}
\item[1.] How long does it take for the procedure to change an arbitrary triangulation of $M$ into a one-vertex
one?

\item[2.] How long does it take for the procedure to find a non-trivial normal $S^{2}$?

\item[3.] How long does it take for the procedure to collapse a non-trivial normal $S^{2}$?

\item[4.] How many non-trivial normal 2-spheres can there be?

\item[5.] How long does it take for the procedure to look for an almost normal $S^{2}$?
\end{enumerate}

\begin{enumerate}
\item[1.] Suppose $M$ has more than one vertex. We take the boundary of a regular neighborhood of an edge. The
weight of that surface is bounded by twice the number of edges. Hence, normalizing this surface takes time linear in
the number of tetrahedra.  We obtain a union of normal 2-spheres. We collapse each of them. See (3.) for the
running of the collapsings.

\item[2.] We have to look through each cone $C_{w, i}(M)$. There are at most $4t \cdot 3^t$ of them. In each
cone, we maximize $\chi$ to find a 2-sphere. This can be done in polynomial time $O(t^n)$, where $n$ is
independent of $t$. Hence the running time of this step is $O(t^{n+1}) \cdot 3^t$.

\item[3.] This part of the procedure refers to Theorem~\ref{MT}. First, to count the number of \R summands from step 2,
 we calculate $H_2(M; \mathbf{Z}) using cell decomposition.$ This can be done in polyniomal time in $t$ by transforming a $(2t) \times t$
  matrix in its row reduced echelon form (see ~\cite{Sch}). It then suffices to compare the number of missing
  $\mathbf{Z}_2$ factors in the second homology of the resulting summands of $M$ with the number of missing $\mathbf{Z}_2$
 factors in the second homology of $M$. Next, we count the number of collapsing surfaces. Note, this number
  depends on the weight of $S$
 but we are only concerned with the surfaces which do not consist entirely of common sides of adjacent $I$-bundles
 or adjacent tips. There cannot be more than 4t such surfaces.
 Finally, it was shown in step 5 that the number of $S^1 \times S^2$ summands is $n + 2  -k$, where $n$ is the number
 of collapsing surfaces which do not consist entirely of common sides of adjacent $I$-bundles or adjacent tips
  and  $k$ is the number of resulting summands. Hence, it takes polynomial time to collapse a
 normal 2-sphere and obtain the triangulations of the resulting pieces as in the conclusion of Theorem~\ref{MT}.

\item[4.] Cutting along a non-trivial normal 2-sphere strictly reduces the number of tetrahedra, so there could
not be more than $t$ of them. In fact, it is shown in ~\cite{Bar3} that there cannot be more than $\lfloor t \rfloor$/2
 of them.

\item[5.] We look through each cone$A_{w, i , l}(M)$. There are at most $3^t \cdot 4t \cdot t$ of them.
In each cone, we maximize $(\chi - o)$. As in 2., such a maximum is found in polynomial time $O(t^m)$,
where $m$ is independent of
$t$. The running time for this step is $O(t^{m+2}) \cdot 3^t$.

\end{enumerate}

\bigskip
\begin{thm} ~\label{C3}
Let M be a closed orientable 3-manifold equipped with a minimal triangulation with t tetrahedra.
 Then there is an algorithm to check if M is reducible or not, and this algorithm runs in polynomial time with
  respect to t.
\end{thm}

\underline{{\bf Proof :}} What makes the above algorithm run in exponential time is the number
of cones $C_{w,i}(M)$ (there are about $3^t$ of them) and $A_{w,i,l}$ (there are about $3^{t}$ of them). When $M$ is
  minimal though, not only we do not need to look through all the cones $C_{w,i}(M)$, but we don't even need to run the
  Thompson-Rubinstein algorithm.

If there are no non-trivial embedded normal 2-spheres, then a famous result of Kneser~\cite{Kn} tells us that
 $M$ does not contain any embedded essential 2-sphere. Hence, $M$ is irreducible. Suppose $M$ does contain a
 non-trivial normal 2-sphere $S$. Let $S$ be a normal 2-sphere constructed in lemma~\ref{L4}.
  When we collapse it we end up with at most 2 summands for $M$.
    Since $M$ is minimal, none of the summands can be homeomorphic to a 3-sphere.
     This tells us that $S$ is essential. Hence, $M$ is reducible. What is important to notice here
     is that $<S> \leq 2$ by Theorem~\ref{MT3}. Therefore, if $M$ is minimal and reducible,
      then it contains an essential normal 2-sphere with 1 or 2 quadrilateral types.

Consider the space of normal surfaces represented by a
family of type $w$, where $w$ assigns the value 0 for the quadrilateral types in every tetrahedra except one. We noted
earlier that surface addition may increase the number of quadrilateral types. Indeed, if $F_{1}$ has exactly one quadrilateral
type in the $j^{th}$ tetrahedra $\Delta_{j}$, and $F_{2}$ has exactly one quadrilateral type in the $i^{th}$
tetrahedra, then $F_{1} + F_{2}$ has a quadrilateral type in $\Delta_{j}$ and $\Delta_{i}$. Hence this space
does not represent a cone. Consider now the space of normal surfaces represented by a family, $w_{(i)}$,
of type $w$,
where $w$ assigns the value 0 for the quadrilateral types in every tetrahedra except in the $i^{th}$ one. It is
now easy to see that this space is a cone in $\mathbf{R}^{7t}$. Every non-trivial normal surface in this cone has one
quadrilateral type, and conversely, every non-trivial normal surface with one quadrilaterl type belong to a cone $C_{w_{(j)},i}(M)$
for some fixed $i$ and $j$. For each
tetrahedron, there are three different quadrilateral types and so there are three $w_{(j)}$'s.
For each type $w_{(j)}$, there are $4t$ choices for the triangle entries $t_{i}$. Hence, there are
$3 \cdot t \cdot (4t) = 12t^{2}$ different cones.

Similarly, we define the space of normal surfaces having two
quadrilateral types in exactly two distinct tetrahedra. If two
tetrahedra $\Delta_j$ and $\Delta_k$ are fixed, this space
represents a cone $C_{w_{(j, k)},i}(M)$. In fact, the union of
these cones, over all possible
 pairs of tetrahedra, represents the space of normal surfaces with one or two quadrilateral types. For each pair of tetrahedra,
  $\{ \Delta_{j}, \Delta_{k} \}$, there are $3 \cdot 3$ possible types $w_{(j, k)}$. There are ${t \choose 2}$
  possible  pairs of distinct tetrahedra. Hence, there are $3 \cdot 3 \cdot {t \choose 2} \cdot (4t) = 36t{t \choose 2}$
   different cones $C_{w_{(j, k)},i}(M)$.

We now describe the algorithm to check if a minimal closed
orientable 3-manifold $M$ is reducible or not. Fix a cone
$C_{w_{(j, k)},i}(M)$ and look at the convex polyhedron $A =
C_{w_{(j, k)},i}(M) \bigcap$ $\{\sum^{7t}_{i=1} t_{i} = 1 \} $. We
maximize $\chi$ on $A$ to obtain a vertex solution S. If $\chi (S)
> 0$, then
 there exists
a 2-sphere in $C_{w_{(j, k)},i}(M)$ and the procedure stops here:
$M$ is reducible. If $\chi (S) \leq 0$, there are no 2-spheres in
the cone $C_{w_{(j, k)},i}(M)$. Repeat this step with a new cone
$C_{w_{(r, s)},i}(M)$. If no 2-spheres have been found in any of
the cones $C_{w_{(j, k)},i}(M)$, then $M$ is irreducible and the
procedure
 stops here. Note, we do not need to look through the cones
 $C_{w_{(j)},i}(M)$ since they lie in some cone $C_{w_{(j, k)},i}(M)$.

As we have seen in Casson's algorithm, it takes polynomial time to
look for
 a 2-sphere in each convex polyhedra $A$. Since there are only $36t{t \choose 2}$ cones of the form $C_{w_{j, k},i}(M)$,
  it will take polynomial time to check if $M$ is reducible or not.

There is, though, a problem in decomposing a minimal 3-manifold in
irreducible pieces in polynomial time. Indeed, after collapsing a
non-trivial normal 2-sphere, we end up with 2 summands which
unfortunately may not be minimal. See ~\cite{Bar2} for a solution
to this problem.

\bigskip

\bigskip

\bibliographystyle{plain}
\bibliography{dd}

\begin{thebibliography}{10}





\bibitem{Bar1} Alexander Barchechat, {\it Minimal Triangulations}, Ph.D. thesis, UC Davis, 2003.

\bibitem{Bar2} Alexander Barchechat, {\it Decomposition of a Closed Orientable Minimal 3-Manifold into irreducible
summands}, in preparation.

\bibitem{Bar3} Alexander Barchechat, {\it A Sharp Upper Bound for the Kneser's Number of Disjoint Essential Non-parallel
2-Spheres in an Closed Orientable 3-Manifold}, in preparation.



\bibitem{Hak} Wolfgang Haken, {\it Theorie der Normalflachen}. Acta Math., 1961, 105, 245-375.

\bibitem{Ha} Joel Hass, {\it Algorithms for Recognizing Knots and 3-Manifolds}. Knot theory and its applications. Chaos Solitons Fractals 9 (1998), no. 4-5, 569--581.

\bibitem{HLP}  Joel Hass, Jeffrey Lagarias and Nicholas Pippenger, {\it The computational Complexity of Knot and Link Problems}. J. ACM 46 (1999), no. 2, 185--211.


\bibitem{JLR} William Jaco, David Letscher, J. Hyam Rubinstein, {\it Algorithms for Essential Surfaces in 3-Manifolds}, in preparation.

\bibitem{JO} William Jaco and Ulrich Oertel, {\it An algorithm to decide if a 3-manifold is a Haken manifold}. Topology 23 (1984), no. 2, 195--209.

\bibitem{JR}  William Jaco and J. Hyam Rubinstein, {\it 0-Efficient Triangulations of 3-Manifolds}. Preprint, 2001.

\bibitem{JS} William Jaco and Eric Sedgwick, {\it Decision Problems in the Space of Dehn Filings}. Topology, 2001.

\bibitem{JT} William Jaco and Jeffrey L. Tollefson, {\it Algorithms for the Complete Decomposition of a Closed 3-Manifold}.  Illinois J. Math. 39 (1995), no. 3, 358--406.

\bibitem{Kar} N. Karmarkar, A new polynomial-time algorithm for linear programming. Combinatorica 4 (1984), no. 4, 373--395.


\bibitem{Kn} H. Kneser,  {\it Geschlossene Flachen in dreidimensionalen Mannifgfaltigkeiten}, Jber. Ent. Math., Ver., 1929, 28, 248-260.

\bibitem{Mat} S. V. Matveev, {\it Classification of sufficiently large three-dimensional manifolds}.  (Russian) Uspekhi Mat. Nauk 52 (1997), no. 5(317), 147--174; translation in Russian Math. Surveys 52 (1997), no. 5, 1029--1055.

\bibitem{Mi}  Aleksandar Mijatovic, {\it Triangulations of Seifert Fibred Manifolds }, arXiv:math.GT/0301246.

\bibitem{Mo}  Edwin E. Moise, {\it Affine Structures in 3-Manifolds: V. The triangulation Theorem and Hauptvermutung}. Ann. of Math. (2) 56, (1952). 96--114.

\bibitem{Pa}  Udo Pachner, {\it Konstruktionsmethoden und das kombinatorische Homomorphieproblem fur Triangulationen kompakter semilinearer Mannigfaltigkeiten}. (German) [Construction methods and the combinatorial homeomorphism problem for triangulations of compact semilinear manifolds] Abh. Math. Sem. Univ. Hamburg 57 (1987), 69--86.


\bibitem{Ru} J. H. Rubinstein, {\it An algorithm to recognize the 3-sphere}.

\bibitem{Sch} Alexander Schrijver, {\it Theory of Linear and Integer Programming}.

\bibitem{Tho} Abigail Thompson, {\it Thin Position and the Recognition Problem for $S^{3}$}. Math. Res. Lett., 1(5):613-630, 1994.










\end{thebibliography}

\end{document}